\documentclass[a4paper,11pt]{article}

\usepackage{theorem}
\usepackage[dvips]{graphics}
\newtheorem{definition}{Definition}
\newtheorem{proposition}{Lemma}
\newtheorem{theorem}{Theorem}

\usepackage{graphics, amssymb, amsmath}
\def\disp{\displaystyle}

\def\mb{\mathbb}
\def\mc{\mathcal}
 
\def\lan{\langle}
\def\ran{\rangle}

\begin{document}
 
\centerline{\textbf{Quadratic BSDEs driven by a continuous martingale}}
\centerline{\textbf{and applications to the utility maximization problem}\footnote{Marie-Amelie Morlais\\ETH Zurich, Switzerland.\\
Ramistrasse 101, HG F27.3, 8006 Zurich,\\ 
Tel: +41 44 632 5859\\ e-mail: marieamelie.morlais@free.fr
}}
\paragraph*{Abstract:} $\quad$
In this paper, we study a class of quadratic Backward Stochastic Differential Equations (BSDEs) which arises naturally in the problem of utility maximization with portfolio constraints. We first establish existence and uniqueness of solutions for such BSDEs and we give applications to the utility maximization problem. Three cases of utility functions: the exponential, power and logarithmic ones, will be discussed.
\section{Introduction}
In this paper, the problem under consideration consists in maximizing the expected utility of the terminal value of a portfolio minus a certain liability and under constraints. The main objective is to give the expression of the value function of the utility maximization problem with utility function $U$ and liability $B$: its expression at time $t$ is
\begin{equation}\label{eq: utilmaxim}
V_{t}^{B}(x) = \textrm{ess} \displaystyle{ \sup_{\nu \in \mc{A}_{t}} \mathbb{E}^{\mathcal{F}_{t}}(U(X_{T}^{\nu, t, x} - B))}.
\end{equation}
In our model, $X_{T}^{\nu, t, x}$ is the terminal value of the wealth process associated with the strategy $\nu$ and equal to $x$ at time $t$ and the $\textrm{ess} \disp{sup}$ is taken over all trading strategies $\nu$
which are in a specific admissibility set $\mc{A}_{t}$ and defined on $[t, T]$. Not any $\mc{F}_{T}$-measurable random variable $B$ is replicable by a strategy taking its values in $\mc{A}_{t}$ and hence, the financial market is incomplete.
This problem provides further interests due to its connection with utility indifference valuation: in fact, the utility indifference price relates the two value processes $V^{B}$ and $V^{0}$. Introduced by Hodges and Neuberger (1989), the utility indifference selling price stands for the amount of money which makes the agent indifferent between selling or not selling the claim $B$.\\ 
Among previous studies of our concerned problem, we can cite \cite{Becherer} and \cite{ManiaetSchw}. In the first, Becherer studies both the utility maximization problem and the notion of utility indifference valuation in a discontinuous setting, whereas, in the second paper, Mania and Schweizer consider the same problem in a continuous framework. As in these two papers and to solve the problem (\ref{eq: utilmaxim}) in the case of non convex trading constraints, we rely on the dynamic programming methodology and on non linear BSDE theory. In the existing literature (see e.g. \cite{Fritelli}, \cite{Owen} or \cite{Schach1}), the convex duality method is widely used to study the unconditional case of the problem, but in the aforementioned papers, the authors either suppose there is no constraints or they assume the convexity of the constraint set, which is an assumption we relax here. We rather use the first method to handle dynamically the problem and, for this approach, the main references are \cite{ImkelleretHu} and \cite{ManiaetSchw}. Our contribution consists in extending the dynamic method in this conti-nuous setting and in presence of constraints. This requires to establish existence and uniqueness results for solutions to specific quadratic BSDEs and then use these results to characterize both the value function expressed in (\ref{eq: utilmaxim}) and the strategies attaining the supremum in the expression of the value function.
\\
The paper is structured as follows: Section 2 lays out the financial background and gives some preliminary tools and results about BSDEs. Then, the dynamic programming method is applied to derive an explicit BSDE. Section 3 investigates the existence and uniqueness results for solutions to the introduced BSDEs. In Section 4, applications to finance are developed and the expression of the value function is provided for three types of utility functions. Lengthy proofs are relegated to an appendix. 

\section{Statement of the problem and main results}
\subsection{The model and preliminaries}
\hspace{0.5cm} We consider ($\Omega$, $\mathbb{F}$, $\mathbb{P}$) a probability space equipped with a right continuous and complete filtration $\mathcal{F}$ = $(\mathcal{F}_{t})_{t}$ and with a continuous $d$-dimensional local martingale $M$.
Throughout this paper, all processes are considered on $[0, T]$, $T$ being a deterministic time and we denote by $Z \cdot M$ the stochastic integral w.r.t. $M$. We assume that $\mathcal{F}$ = $(\mc{F}_{t})_{t \in [0, T]}$ is a continuous filtration: this means that any martingale $K$ is continuous. Any $\mathbb{R}$-valued square integrable martingale $K$ can be written 
$$K = Z \cdot M + L,$$ 
with $Z$ a predictable $\mathbb{R}^{d}$-valued process and $L$ a (square integrable) $\mb{R}$-valued martingale strongly orthogonal to $M$ (i.e., for each $i$, $\lan{M}^{i}, L \ran$ = 0).
 For a given square integrable martingale $M$, the notation $\lan{M \ran}$ stands for the quadratic variation process.\\
\indent  From the Galchouk-Kunita-Watanabe inequality, it follows that each component $d\lan{M}^{i},M^{j}\ran$ ($i, j \in \{1, \cdots , d\}$) is absolutely continuous with respect to $d\tilde{C} = (\mathop{\sum_{i}d\lan{M}^{i}\ran })$. 
Hence, there exists an increasing and bounded process $C$ (for instance, we set: $C_{t}$ = $\arctan(\tilde{C}_{t})$, for all $t$) such that $\lan{ M \ran} $ can be written
$$ d\lan{M} \ran _{s} = m_{s} m_{s}^{'}dC_{s},$$
where $m$ is a predictable process taking its values in $\mathbb{R}^{d \times d}$ (this expression has been used in \cite{ElkaretMazliak} in an analogous continuous framework). The notation $m^{'}$ stands for the transposed matrix and we also assume that, for any $s$, the matrix $m_{s}m_{s}^{'}$ is invertible.
\paragraph*{The financial background}\hspace{0.5cm}
 To bring further motivations, we explain the financial context and, for this, we provide here all the definitions and common assumptions.
We consider a financial market consisting in $d + 1$ assets: one risk free asset with zero interest rate and $d$ risky assets. We model the price process $S$ of the $d$ risky assets as a process satisfying
\begin{equation}\label{eq: SDE1}
 \frac{dS_{s}}{S_{s}} = dM_{s} + dA_{s}, \; \; \textrm{with:} \;\forall \; j \in \{1, \cdots,d\},\; dA^{j}_{s} = \disp{\sum_{i=1}^{d}\lambda_{s}^{i}d\lan{M^{j}}, M^{i} \ran_{s}}, 
\end{equation}
 $(\lan{M^{j}},M^{i} \ran)_{j}$ standing for the $i^{th}$ column of the $\mathbb{R}^{d \times d}$ matrix-valued process $\lan{M} \ran$.
This decomposition already introduced in \cite{DelbaenetSch} and the assumption of the almost sure inversibility of $m_{s}$ for all $s$,  
 ensure that the usual no arbitrage property holds. We also assume that $\lambda$ is a $\mathbb{R}^{d}$-valued process satisfying
\begin{equation}\label{eq: hyplambda}
(H_{\lambda}) \quad  \quad \quad \quad  \exists \; a_{\lambda} > 0, \; \; \int_{0}^{T}\lambda_{s}^{'}d\lan{M} \ran_{s}\lambda_{s} = \int_{0}^{T}|m_{s}\lambda_{s}|^{2}dC_{s} \leq a_{\lambda}, \quad \mb{P}\textrm{-a.s.}   
\end{equation}
This definition is stronger than the usual structure condition, which only states:\\ $\disp{\int_{0}^{T}\lambda_{s}^{'}d\lan{M} \ran_{s}\lambda_{s}} < \infty,$ $\mb{P}$-a.s. (we refer to \cite{Stricker2} or \cite{HPSCHW01} for this condition) and implies that $\mc{E}(- \lambda \cdot M)$ is a strict martingale density for the price process $S$. In the financial application, we rely on ($H_{\lambda}$) to use the precise a priori estimates given in Lemma \ref{estimapriori} in Section 3. We now state the definition of wealth process $X^{\nu}$ and of the associated self-financing and constrained trading strategy $\nu$.
\begin{definition}
A predictable $\mathbb{R}^{d}$-valued process $\nu:= (\nu_{s})_{s \in [t, T]}$
 is called a self-financing trading strategy if it satisfies
 \begin{itemize}
  \item[1]$\nu_{s} \in \mc{C}$, $\mb{P}$-a.s. and for all $s$, $\mc{C}$ being the constraint set (closed and not necessarily convex set in $\mb{R}^{d}$).\\
 \item[2]The wealth process $X^{\nu}: = X^{\nu, t, x}$ of an agent
with strategy $\nu$ and wealth $x$ at time $t$ is defined as \\
 \begin{equation}\label{eq: wealthproc}
 \disp{ \forall s \in \left[t ,T\right], \quad  X_{s}^{\nu} = x + \int_{t}^{s}\mathop{\sum_{i=1}^{d}\frac{\nu_{r}^{i}}{S_{r}^{i}}dS_{r}^{i}},}
  \end{equation}
 and it is in the space $\mc{H}^{2}$ of semimartingales.
\end{itemize}
 \end{definition}

\indent In this definition, each component $\nu_{i}$ of the trading strategy corresponds to the amount of money invested in the $i^{th}$ asset.
Due to the presence of portfolio constraints, there does not necessarily exist a strategy $\nu$ (such that, for all $s$, $\nu_{s} \in \mc{C}$) and satisfying:
$X_{T}^{\nu} = B,$ for a given $\mc{F}_{T}$-measurable claim. Hence, we are facing an incomplete market.
 The utility maximization problem aims at giving the expression of the value function defined in (\ref{eq: utilmaxim}) and at characterizing the set of optimal strategies, i.e. those achieving the $\textrm{ess} \disp{\sup}$ for the problem. In our study, we first consider the exponential utility maximization problem associated with the utility function:
$ U_{\alpha}(x) = -\exp{(-\alpha x)}.$
Usually, the set of admissible trading strategies consists of all the strategies such that the wealth process is bounded from below. 
To solve the problem, we need to enlarge the set of admissible strategies analogously to \cite{ImkelleretHu} to a new set denoted by $\mc{A}_{t}$.
\begin{definition}
Let $\mathcal{C}$ be the constraint set, which is such that: $0 \in \mc{C}$. The set $\mathcal{A}_{t}$ of admissible strategies consists of all $d$-dimensional predictable processes $\nu := (\nu_{s})_{s \in [t, T]}$ satisfying $\nu_{s} \in \mathcal{C}$, $\mb{P}$-a.s. and for all $s$, $\mathbb{E}(\disp{\int_{t}^{T}|m_{s}\nu_{s}|^{2}dC_{s}}) < \infty $ and the uniform integrability of the family\\
$$ \{ \exp(-\alpha X_{\tau}^{\nu}) : \;  \tau \; \mc{F}\textrm{-stopping time taking its values in}\; [t, T] \} .$$
\end{definition}
 This appears to be a restrictive condition on strategies. Hence, we have to justify the existence of one optimal strategy admissible in this sense.

\paragraph*{Preliminaries on quadratic BSDEs}\hspace{0.5cm}
 We first provide the form of the one dimensional BSDEs considered in the sequel 
 \[ (\textrm{Eq1}) \quad  
Y_{t} = B +  \int_{t}^{T}F(s,Y_{s}, Z_{s})dC_{s}  + \frac{\beta}{2}(\lan{L}\ran_{T} - \lan{L}\ran_{t}) -
\int_{t}^{T}Z_{s}dM_{s} - \int_{t}^{T}dL_{s}.    \]
To refer to this BSDE, we use the notation BSDE$(F, \beta, B)$. Usually, a BSDE is cha-racterized by two parameters: its terminal condition $B$ assumed here to be bounded, its generator $F:=F(s, y, z)$, a $\mc{P} \times \mc{B}(\mb{R}) \times \mc{B}(\mb{R}^{d})$-measurable function, continuous w.r.t. ($y, z$) ($\mathcal{P}$ denotes the $\sigma$-field of all predictable sets of $[0, T] \times \Omega$ and $\mathcal{B}(\mathbb{R})$ the Borel field of $\mb{R}$). In our setting, we introduce another parameter $\beta$ which is assumed to be constant and a financial meaning for $\beta$ is given in next paragraph. We will also impose precise growth conditions on the generator (and, in particular, the quadratic growth w.r.t. $z$). One essential motivation of this study is that such quadratic BSDEs \footnote{Such BSDEs have been considered in \cite{ManiaetSchw}, where the authors deal with the utility maximization problem but, contrary to the present paper, they assume there is no trading constraints.} appear naturally when using the same dynamic method as in \cite{ImkelleretHu} to solve the problem (\ref{eq: utilmaxim}). 
 A solution of the BSDE$(F, \beta, B)$ is a triple of processes ($Y, Z, L$) with $\lan{L}, M \ran = 0$, such that: $\disp{\int_{0}^{T} |F(s, Y_{s}, Z_{s})|dC_{s}} < \infty, \; \mb{P}$-a.s., satisfying ($\textrm{Eq1}$) and defined on $S^{\infty} \times L^{2}(d \lan M \ran \otimes d\mathbb{P}) \times \mathcal{M}^{2}([0, T])$:
  $S^{\infty}$ consists of all bounded continuous processes, $L^{2}(d \lan M \ran \otimes d\mathbb{P})$ consists of all predictable processes $Z$ such that:
$\disp{\mathbb{E}\big(\int_{0}^{T}|m_{s}Z_{s}|^{2}dC_{s}\big) < \infty}$, and $\mathcal{M}^{2}([0, T])$
consists of all real square integrable martingales of the filtration $\mathcal{F}$.\\
\indent 
  The stochastic exponential denoted by $\mc{E}(K)$ of a semimartingale $K$ is the unique process satisfying
$$\mathcal{E}_{t}(K) = 1 + \int_{0}^{t}\mathcal{E}_{s}(K)dK_{s}.$$
 A process $L$ is a BMO martingale if $L$ is a $\mathcal{F} $ martingale and if there exists a constant $c$ ($c > 0$) such that, for any $\mathcal{F}$ stopping time $\tau $,
$$ \quad
\mathbb{E}^{\mathcal{F}_{\tau}}(\lan{L} \ran _{T} - \lan{L} \ran _{\tau} ) \leq c . $$

\paragraph*{The dynamic method}\hspace{0.5cm}
 In this part, we use the same dynamic method as in \cite{ImkelleretHu} to characterize the value function of the optimization problem in terms of the solution of a BSDE with parameters ($F^{\alpha}$, $\beta$, $B$). The expressions of $F^{\alpha}$ and $\beta$ are obtained below by formal computations (these computations are justified in the last section of this paper).\\
\indent For this, we construct,
for any strategy ${\nu}$ and fixed $t$, a process $R^{\nu}:= (R_{s}^{\nu})_{s \ge t}$ such that, for all $s$,
  $R_{s}^{\nu} = U_{\alpha}(X_{s}^{\nu} - Y_{s})$, with $U_{\alpha}$ defined by: $U_{\alpha}(\cdot): = -\exp(-\alpha \cdot)$), and where the process $Y$ is a solution of a BSDE($F^{\alpha}$, $\beta$, $B$) of type ($\textrm{Eq1}$): the terminal condition is the contingent claim $B$, and the parameters $F^{\alpha}$ and $\beta$ have to be determined.
Besides, this family ($R^{\nu}$) is such that
\\
(i) $R_{T}^{\nu}$ = $U_{\alpha}(X_{T}^{\nu} - B)$, for any strategy $\nu$,\\
(ii) $R_{t}^{\nu}$ = $R_{t}$ = $U_{\alpha}(x - Y_{t})$ ($x$ is assumed to be a constant\footnote{This dynamic method can be extended to any attainable wealth $x$, i.e. any $\mc{F}_{t}$-measurable random variable such that: $X_{t}^{\pi} :=x$ for at least one admissible strategy $\pi$ defined on $[0,t]$.}).\\
(iii) $R^{\nu}$ is a supermartingale for any strategy $\nu$, $\nu \in \mc{A}_{t}$, and a martingale for a parti-cular strategy $\nu^{*}$, $\nu^{*} \in \mc{A}_{t}$.\\
We rely on the equation (\ref{eq: wealthproc}) defining $X^{\nu}$ and on It\^o's formula to get
\begin{tabbing}
 $\disp{X_{s}^{\nu} - Y_{s} =}$ \=$\disp{ (x - Y_{t}) + \int_{t}^{s}(\nu_{u} - Z_{u})dM_{u} - (L_{s} -L_{t}) }$ \\
\> $\disp{+ \int_{t}^{s}F^{\alpha}(u, Z_{u})dC_{u} + \frac{\beta}{2}(\lan{L} \ran_{s} -  \lan{L} \ran_{t})+ \int_{t}^{s}(m_{u}\nu_{u})^{'}(m_{u}\lambda_{u})dC_{u}}.$ \\
\end{tabbing}
Since, for all $s$, $R_{s}^{\nu}: = - \exp(-\alpha(X_{s}^{\nu} - Y_{s}))$ and using the notation: $\mc{E}_{t, T}(K) = \frac{\mc{E}_{T}(K)}{\mc{E}_{t}(K)},$ for a given local martingale $K$, we claim
\[  \begin{array}{ll}
 \exp(- \alpha(\disp{\int_{t}^{T}(\nu_{s}- Z_{s})dM_{s})}) =  \quad \quad &\\ 
   \mathcal{E}_{t, T}(- \alpha((\nu- Z) \cdot M))\exp(\frac{\alpha^{2}}{2}\disp{\int_{t}^{T}|m_{s}(\nu_{s} - Z_{s})|^{2}dC_{s}}), & \textrm{ on the one hand,} \\
\\ 
 \disp{\exp( \alpha (L_{T} - L_{t})) = \mathcal{E}_{t, T}(\alpha L)\exp(\frac{\alpha^{2}}{2}(\lan{L} \ran_{T} - \lan{L} \ran_{t}) ), } \quad  & \textrm{ on the other hand,}\\
\end{array} \]
  which leads to the mutiplicative decomposition 
\begin{equation}\label{eq: decompmultip}  R_{s}^{\nu} = -\exp(-\alpha(x- Y_{t}))\mathcal{E}_{t, s}(-\alpha(\nu - Z) \cdot M)\mathcal{E}_{t, s}(\alpha L)\exp(A_{s}^{\nu} - A_{t}^{\nu}).\end{equation}
Here, $A^{\nu}$ is such that
\[ \begin{array}{lc}  \disp{dA_{s}^{\nu} =} & \disp{(-\alpha F^{\alpha}(s, Z_{s}) - \alpha(m_{s}\nu_{s})^{'}\big(m_{s}\lambda_{s}) +\frac{\alpha^{2}}{2}|m_{s}(\nu_{s} - Z_{s})|^{2}\big)dC_{s} } \\   
&\disp{ + (\frac{\alpha^{2} - \alpha \beta}{2})d\lan{L} \ran_{s}} \; . \\
\end{array} \]
$M$ and $L$ being strongly orthogonal, $$\disp{ \mathcal{E}(-\alpha(\nu - Z) \cdot M)\mathcal{E}(\alpha L) = \mathcal{E}(-\alpha(\nu - Z) \cdot M + \alpha L)}. $$
In (\ref{eq: decompmultip}), $R^{\nu}$ is the product of a positive local martingale (as a continuous stochastic exponential) and a finite variation process. $R_{t}^{\nu}$ being negative $\mb{P}$-a.s. and relying on the multiplicative decomposition (\ref{eq: decompmultip}), the increasing property of $A^{\nu}$ for all $\nu$ yields the supermartingale property of $R^{\nu}$ (it is a martingale for $\nu^{*}$ satisfying: $dA^{\nu^{*}} \equiv 0$).
These two last conditions on the family $(A^{\nu})$ holding true for all $\nu$, $\nu \in \mc{A}_{t}$, we get
\[ \left\{ \begin{array}{l} -\alpha \frac{\beta}{2}d\lan{L} \ran_{s} + \frac{\alpha^{2}}{2}d\lan{L} \ran_{s} = 0 \quad  \Rightarrow\;  (\beta = \alpha), \;  \\
 - \alpha(F^{\alpha}(s, Z_{s})+ (m_{s}\nu_{s})^{'}(m_{s}\lambda_{s}) ) + \frac{\alpha^{2}}{2}|m_{s}(\nu_{s}-Z_{s})|^{2} \ge 0.   \end{array} \right.  \]
This leads to
\begin{eqnarray}\label{eq: parameterF}
 F^{\alpha}(s, z) = \displaystyle{  \mathop{\inf_{\nu \in \mathcal{C}}\big(\frac{\alpha}{2}|m_{s}(\nu - (z +  \frac{\lambda_{s}}{\alpha}))|^{2}\big) }}
- (m_{s}z)^{'}(m_{s}\lambda_{s})  - \frac{1}{2\alpha}|m_{s}\lambda_{s}|^{2}.  \\
\nonumber  \end{eqnarray}
This method, explained for a fixed time $t$, relies on the dynamic programming principle and could therefore be extended without any additional difficulty to any $\mc{F}$-stopping time $\tau$.
\subsection{Statement of the assumptions and main results}
\paragraph*{Assumptions}\hspace{0.5cm}
To study the existence for solutions of the BSDEs($F, \beta , B$) of type ($\textrm{Eq1}$), we assume in all the sequel the boundedness of the terminal condition $B$. Moreover, we use one of the three following assumptions:
we suppose that there exists a predictable process denoted by $\bar{\alpha} $ such that
\begin{equation} 
 \;\bar{\alpha} \geq 0 \; \textrm{and}\; \int_{0}^{T}\bar{\alpha}_{s}dC_{s} \leq a, \; (a > 0),\;  \exists \; b, \; \gamma , \; C_{1} > 0,   
\end{equation} 
\[ \begin{array}{lc}
 (H_{1})   & \disp{|F(s, y, z)| \leq \bar{\alpha}_{s}+  b\bar{\alpha}_{s}|y| + \frac{\gamma}{2}|m_{s}z|^{2}} \; \disp{\textrm{with} \; \gamma \geq |\beta| \; \textrm{and } \;  \gamma \geq b, } \\

  \\
\\
 (H_{1}^{'}) &  \disp{  |F(s, y, z)| \leq  \bar{\alpha}_{s} + \frac{\gamma}{2}|m_{s}z|^{2}}, \\
\\
 (H_{1}^{''}) & \disp{- C_{1}( \bar{\alpha}_{s} + |m_{s}z| )   \leq  F(s, y, z) \leq    \bar{\alpha}_{s} + \frac{\gamma}{2}|m_{s}z|^{2}}. \\
\end{array} \] 
\textbf{Remark:}
$\bullet\;$Assumption $(H_{1})$ is more general than the two other ones but we will only require these two last assumptions to establish the existence result. We first reduce the assumption ($H_{1}$) to ($H_{1}^{'}$) by a classical truncation procedure and we note that the additional assumption in $(H_{1}^{''})$ is that the lower bound has at most linear growth in $z$: this condition has already been used by \cite{BriandetHu} in the Brownian setting to justify the existence of a minimal solution. We rely on the same construction to prove our existence result.\\
 $\bullet\;$ The quadratic BSDE introduced in Section 2.1 and of the form ($\textrm{Eq1}$) has for parameters $F:= F^{\alpha}$, $\beta :=\alpha$ and for terminal condition the contingent claim $B$. In particular, the generator $F^{\alpha}$ given by (\ref{eq: parameterF}) satisfies ($H_{1}$).
In fact, we have
\begin{eqnarray*}
F^{\alpha}(s, z) \geq - (m_{s}z)^{'}(m_{s}\lambda_{s})  - \frac{1}{2\alpha}|m_{s}\lambda_{s}|^{2}
\geq -  |m_{s}z||m_{s}\lambda_{s}| - \frac{1}{2\alpha}|m_{s}\lambda_{s}|^{2} , \\
 \end{eqnarray*}
which leads to
\begin{eqnarray*}
F^{\alpha}(s, z) \geq - \big(\frac{\alpha}{2}|m_{s}z|^{2} + \frac{1}{\alpha}|m_{s}\lambda_{s}|^{2}\big).
\end{eqnarray*}
Defining $\bar{\alpha}$, for all $s$, by: $\bar{\alpha}_{s} = \frac{1}{\alpha}|m_{s}\lambda_{s}|^{2}$, we claim that: $\disp{\int_{0}^{T}\bar{\alpha}_{s}dC_{s}} \leq a$, $\mb{P}$-a.s., with the parameter $a$ depending on $\alpha$ and $a_{\lambda}$ (defined in equation (\ref{eq: hyplambda}) referred as ($H_{\lambda}$)). Noting that 0 is in $ \mathcal{C}$, we get
\[ F^{\alpha}(s, z) \leq \frac{\alpha}{2}|m_{s}z|^{2} .\]

$\bullet\;$ Even if we suppose that $F$ is Lipschitz w.r.t. $y$ and $z$, we cannot obtain directly existence and uniqueness result for a BSDE of type ($\textrm{Eq1}$), because of the presence of the additional term involving the quadratic variation process $\lan{L} \ran$. This explains the
introduction of another type of BSDEs denoted by ($\textrm{Eq2}$)
\[ (\textrm{Eq2}) \left\{ \begin{array}{lc}
dU_{s} = - g(s, U_{s}, V_{s})dC_{s} + V_{s}dM_{s} + dN_{s},  \\
U_{T} = e^{\beta B}. \quad \quad \quad \quad &\\  
\end{array}  \right. \]
In the sequel, we denote it by BSDE($g, e^{\beta B}$). This second type of BSDE is linked with the BSDE($F, \beta, B$) of type ($\textrm{Eq1}$) by using an exponential change of variable. Hence, setting: $U:= e^{\beta Y}$, this leads to
$$
 g(s, u, v) = \big(\beta u F(s, \frac{\ln(u)}{ \beta},  \frac{v}{\beta u}) -  \frac{1}{2u}|m_{s}v|^{2}\big) \mathbf{1}_{u > 0}.$$ 
This second type of BSDE is simpler, since there is no more term involving the quadratic variation process $\lan{N} \ran $ in $(\textrm{Eq2})$. This type of BSDEs with $g$ uniformly Lipschitz w.r.t. $y$ and $z$ has been studied in \cite{ElkaretMazliak} in a general continuous setting. Our aim is to establish a one to one correspondence between the solutions of the BSDE($F, \beta, B$) of type ($\textrm{Eq1}$) and those of the BSDE($g, e^{\beta B}$) of type ($\textrm{Eq2} $).\\

\indent To prove a uniqueness result for solutions of the BSDE($F, \beta, B$) of type ($\textrm{Eq1}$), we impose that there exists two reals $\mu$ and $C_{2}$, a non negative predictable process $\theta$ and a constant $c_{\theta}$ such that 
\[(H_{2}) \left\{ \begin{array}{lr}
   \forall \; z \in \mathbb{R}^{d},\; \forall \;y^{1}, \; y^{2} \in \mathbb{R},\;   \\
\\
(y^{1} - y^{2})(F(s, y^{1}, z) - F(s, y^{2}, z)) \leq \mu|y^{1} - y^{2}|^{2}, \\
\\
  \exists \; \theta \;  \textrm{s.t.}\; \disp{\int_{0}^{T}|m_{s}\theta_{s}|^{2}dC_{s}} \leq c_{\theta}, \;  \forall \; y \in \mathbb{R}, \;\forall \;  z^{1},\; z^{2} \in \mathbb{R}^{d},&     \\
\\
|F(s, y, z^{1}) - F(s, y, z^{2})| \leq C_{2}(m_{s}\theta_{s} + |m_{s}z^{1}| + |m_{s}z^{2}|)|m_{s}(z^{1} - z^{2})|. \\

\end{array} \right. \]
\textbf{Remark:}
The first inequality in assumption ($H_{2}$) corresponds to the monotonicity assumption (this assumption is given in \cite{Pardoux2}).
The second assumption on the increments in the variable $z$ is a kind of local Lipschitz condition w.r.t $z$, which is similar to the one in \cite{ImkelleretHu}.
 We check that ($H_{2}$) is satisfied by the generator $F^{\alpha}$ with: $C_{2}= \frac{\alpha}{2}$, $\theta \equiv 4\frac{|m \lambda|}{\alpha }$ and $\mu = 0$, since $F^{\alpha} $ is independent of $y$. For any $z^{1}$, $z^{2}$ in $\mathbb{R}^{d}$, the expression of the increments in $z$ is
\[ \begin{array}{ll}
\disp{|F^{\alpha}(s, z^{1})} & \disp{-\;F^{\alpha}(s, z^{2})| } \\
\\ 
& \disp{ \leq |\frac{\alpha}{2}\big(\mathop{\textrm{dist}^{2}(m_{s}(z^{1} + \frac{\lambda}{\alpha}), m_{s}\mathcal{C})} - \mathop{\textrm{dist}^{2}(m_{s}(z^{2} + \frac{\lambda}{\alpha}), m_{s}\mathcal{C})}\big)|}  \\
&\quad  \disp{ \; \; + |-(m_{s}z^{1})^{'}(m_{s}\lambda) + (m_{s}z^{2})^{'}(m_{s}\lambda)| } \\
\\
 & \disp{\leq \frac{\alpha}{2}|m_{s}(z^{1} - z^{2})| \big(|m_{s}z^{1}| + |m_{s}z^{2}| + 2\frac{|m_{s}\lambda|}{\alpha}\big) +  |m_{s}(z^{1} - z^{2})||m_{s}\lambda|}.\\
\end{array} \]
 
\paragraph*{ Main results}
\hspace{0.5cm} To obtain the existence and uniqueness results for solutions of BSDEs of type ($\textrm{Eq1}$), we  establish the same results for BSDEs of type ($\textrm{Eq2}$). Note that the BSDEs of type ($\textrm{Eq2}$) are a particular case of the BSDEs of the first type without the additional term in $\lan{ L} \ran$. We now state the results which are justified in Section 3.
\begin{theorem}\label{existence}\textbf{Existence:}
 Considering the BSDE$(F, \beta, B)$ and assuming that the generator $F$ satisfies ($H_{1}$), there exists a solution ($Y$, $Z$, $L$) in $S^{\infty} \times L^{2}(d\lan{M} \ran \otimes d\mathbb{P}) \times \mathcal{M}^{2}([0,T])$ of the BSDE.\\ 
\end{theorem}


\begin{theorem}\label{unicite}\textbf{Uniqueness:}
For all BSDEs($F, \beta, B$) of type ($\textrm{Eq1}$) such that the generator $F$ satisfies both ($H_{1}$) and ($H_{2}$) and such that the terminal condition is bounded, there exists a unique solution ($Y$, $Z$, $L$) in $S^{\infty} \times L^{2}(d\lan{M} \ran \otimes d\mathbb{P}) \times \mathcal{M}^{2}([0,T])$.
\end{theorem}

\begin{theorem}\label{comparison}\textbf{Comparison:}
Considering two BSDEs of the form $(\textrm{Eq1})$ given by ($F^{1}, \beta, \xi^{1}$) and ($F^{2}, \beta, \xi^{2}$) and such that $F^{1}$ and $F^{2}$ satisfy ($H_{1}$) and ($H_{2}$) and assuming furthermore that ($Y^{1}$, $Z^{1}$,
$L^{1}$) and ($Y^{2}$, $Z^{2}$, $L^{2}$) are respective solutions of each BSDE such that 
  $$\left( \xi^{1} \leq \xi^{2}  \quad \textrm{and} \quad
 F^{1}(s, Y_{s}^{1}, Z_{s}^{1}) \leq F^{2}(s, Y_{s}^{1}, Z_{s}^{1})\right), \; \;\mb{P}\textrm{-a.s. and for all} \;s,$$
then, we have: $\; Y_{s}^{1} \leq Y_{s}^{2}$, $\mb{P}$-a.s. and for all $s$.\\
\end{theorem}
We only provide proofs for the two first theorems, since,
without additional difficulty, we check that the comparison result given in Theorem \ref{comparison} holds: to prove it, we proceed with a linearization of the generator similar as the one in Section 3.2: this consists in applying It\^o-Tanaka formula to the process: $\tilde{Y}_{\cdot}^{1, 2} := \exp(2\mu C_{\cdot})|(Y_{\cdot}^{1} - Y_{\cdot}^{2})^{+}|^{2}$ and then in rewriting identically the proof.

\section{Results about quadratic BSDEs }

\subsection{A priori estimates }
\hspace{0.5cm}In this part, we obtain precise a priori estimates for solutions of the BSDEs of type ($\textrm{Eq1}$). Referring to previous studies on quadratic BSDEs (such as in \cite{BriandetHu} or \cite{mkobylanski}), these estimates are the starting point of the proof of the main existence result.\\
 To prove these estimates, we assume the existence of a solution ($Y$, $Z$, $L$) of the BSDE($F, \beta, B$) such that $F$ satisfies ($H_{1}$)) and we proceed analogously to \cite{BriandetHu}. However, since the authors work with a brownian filtration, we have to generalize their method to our setting.
\begin{proposition}\label{estimapriori}
Considering a BSDE of type ($\textrm{Eq1}$) given by ($F, \beta, B$) and assuming both boundedness of $B$ and condition ($H_{1}$) for $F$, there exists three positive constants $ c,\; C, \; C^{'}$ depending only on $\gamma$, $a$, $b$ (given in ($H_{1}$)) and on $|B|_{\infty}$ such that, for any solution ($Y, Z, L$), 
  \[ 
\begin{array}{l} 
  (i) \; \disp{ \mathbb{P}\textrm{-a.s. and for all}\;
             t, \; c \leq  Y_{t} \leq C, }  \\
            
   (ii)  \; \disp{ \textrm{for any} \; \mathcal{F}
             \textrm{-stopping time} \; \tau, \;
      \mathbb{E}^{\mathcal{F}_{\tau}}\left(\int_{\tau}^{T}|m_{s}Z_{s}|^{2} dC_{s}+ \lan{L} \ran_{T}  - \lan{L}\ran_{\tau} \right)
             \leq  C^{'}}. \\
\end{array} \]
 \end{proposition}
 By definition, the solution ($Y, Z, L$) belongs to $S^{\infty} \times L^{2}(d \lan{M} \ran \otimes d\mathbb{P}) \times \mathcal{M}^{2}([0,T])$: in particular, $Y_{s}$ is almost surely bounded (uniformly in $s$, $s$ in $[0, T]$).
Now, for any $z$, $ z \ge 0$, we introduce $\phi_{t}(z)$ which is defined on $[0, T]$ as 
\begin{equation}\label{eq: fonctionphi} \disp{\forall \;t,\quad \phi_{t}(z) = \exp(\gamma\frac{e^{\int_{t}^{T}b\bar{\alpha}_{u}dC_{u}} - 1}{b})\exp(\gamma z e^{\int_{t}^{T}b\bar{\alpha}_{u}dC_{u}}),} 
\end{equation}
and we check
$$ \forall \; t \in [0, T], \quad 1 \le \phi_{T}(z) \le \phi_{t}(z) \le \phi_{0}(z).$$
Setting: $ \tilde{a} := \frac{e^{ba} -1 }{b}$, we aim at proving 
\begin{equation}\label{eq: domination}\forall \; t, \quad  \exp(\gamma |Y_{t}|) \le \mb{E}\big(\phi_{0}(|B|) | \mc{F}_{t}\big) \le \exp\big(\gamma(\tilde{a} + |B|_{\infty}e^{ba}) \big).\end{equation}
For this, we introduce $H := \big(U(s, |Y_{s}|)\big) $ such that
 $$H_{s} = U(s,|Y_{s}|) := \exp\left( \gamma \big( \frac{\exp(\int_{t}^{s}b \bar{\alpha}_{u}dC_{u}) -1}{b}\big) + \gamma |Y_{s}| \exp(\int_{t}^{s}b \bar{\alpha}_{u}dC_{u})\right), $$
and satisfying: $ U(t, |Y_{t}|) = e^{\gamma |Y_{t}|}$. 
Applying It\^o's formula to the process $H$, we justify that it is a local submartingale: to this end, we prove that the predictable bounded va-riation process $A$ in the canonical decomposition of the semimartingale $H$ is increasing.
 For sake of clarity, we first apply the It\^o-Tanaka formula to $|Y|$
\[ \begin{array}{ll}
d|Y_{s}| = & -\textrm{sign}(Y_{s})F(s, Y_{s}, Z_{s})dC_{s} -\textrm{sign}(Y_{s})\frac{\beta}{2}d\lan{L} \ran_{s} + d\ell_{s}\\
\\
& \; + \textrm{sign}(Y_{s})\big( Z_{s}dM_{s} + dL_{s}\big) , \\
\end{array} \]
 $\ell$ being the local time of $Y$.
Now, It\^o's formula leads to the following expression of $A$
\[ \begin{array}{l} \exp(-\int_{t}^{s}b \bar{\alpha}_{u}dC_{u})dA_{s}  :=  \\
\\
\quad \quad \quad  H_{s}\left( \gamma \bar{\alpha}_{s} - \gamma \textrm{sign}(Y_{s})F(s, Y_{s}, Z_{s}) + \gamma b \bar{\alpha}_{s} |Y_{s}| + \frac{\gamma^{2}}{2}e^{\int_{t}^{s}b \bar{\alpha}_{u}dC_{u}}|m_{s}Z_{s}|^{2}  \right)dC_{s}  \\
    \\
\quad \quad \quad + \;H_{s}\gamma d\ell_{s} 
+\; H_{s}\gamma \left(\big(\frac{\gamma}{2}\exp(\int_{t}^{s}b \bar{\alpha}_{u}dC_{u}) - \textrm{sign}(Y_{s})\frac{\beta}{2} \big)d\lan{L} \ran_{s} \right). \\
\end{array} \]
Using assumption ($H_{1}$) and the inequalities: $|\beta| \le \alpha$ and $\bar{\alpha} \ge 0$, we get that:
$\big(A_{t} := \disp{\int_{0}^{t}dA_{s}} \big)$ is an increasing process. Hence, $H$ is a local submartingale and we conclude relying on a standard localization procedure: 
there exists a sequence $(\tau_{k})$ of increasing stopping times, converging to $T$ and taking values in $[t, T]$ and such that $ (U(s \wedge \tau_{k}, |Y_{s \wedge \tau_{k}}|) $ is a submartingale. This entails 
$$ \quad e^{\gamma |Y_{t}|} = U(t, |Y_{t}|) \leq \mb{E}\big( U(T\wedge \tau_{k}, |Y_{T \wedge \tau_{k}}|)| \mc{F}_{t}\big). $$ 
Applying the bounded convergence theorem to $\left(\mb{E}\big( U(T\wedge \tau_{k}, |Y_{T \wedge \tau_{k}}|)| \mc{F}_{t}\big)\right)_{k}$ 
and letting $k$ tend to infinity, we obtain
$$e^{\gamma |Y_{t}|} \le \mb{E}\big(U(T, |Y_{T}|) | \mc{F}_{t}\big), $$
which gives (\ref{eq: domination}).
 Hence, assertion (i) of lemma \ref{estimapriori} is satisfied with
$$ C := (\tilde{a} + |B|_{\infty}e^{ba})\; \textrm {and} \; c  := -(\tilde{a} + |B|_{\infty}e^{ba}).$$

\indent To prove assertion (ii), we apply It\^o's formula to the bounded process $\tilde{\psi}(Y) :=  \psi_{\gamma}(Y + |c|)$, with $\psi_{\gamma}$ such that
$$ \psi_{\gamma}(x) = \frac{e^{\gamma x} - 1 - \gamma x}{\gamma^{2}}.$$
$c$ being the lower bound of $Y$, we have: $Y + |c| \geq 0,$ $\mb{P}$-a.s. We also use the properties
  \begin{equation}\label{eq: Proppsi}  \psi_{\gamma}{'}(x) \geq 0,  \;  \textrm{if}\; x \geq 0, \; \; \textrm{and} \; \;
-\gamma \psi_{\gamma}^{'} + \psi_{\gamma}^{''} = 1.  \\
 \quad \quad \end{equation}
We now consider an arbitrary stopping time $\tau$ of ($\mathcal{F}_{t})_{t \in [0, T]}$. Taking the conditional expectation with respect to $\mathcal{F}_{\tau }$ in It\^o's formula given between $\tau$ and $T$, we get\\
\begin{tabbing}
 $ $  \=$\disp{\tilde{\psi}(Y_{\tau}) -} \disp{ \mathbb{E}^{\mathcal{F}_{\tau}}(\tilde{\psi}(Y_{T})) }$ \\
  $ \quad \;$ \> $
   \disp{ = -\mathbb{E}^{\mathcal{F}_{\tau  }}\left(\int_{\tau }^{T }\tilde{\psi}^{'}(Y_{s})
  ( - F(s,  Y_{s}, Z_{s})dC_{s} - \frac{\beta}{2}d\lan{L} \ran_{s}) \right)} $\\ 
\> $ \quad \disp{-\; \mathbb{E}^{\mathcal{F}_{\tau }}\left(\int_{\tau }^{T}\tilde{\psi}^{'}(Y_{s })(Z_{s}dM_{s} + dL_{s})\right) } \disp{-\; \mathbb{E}^{\mathcal{F}_{\tau }}\left( \int_{\tau }^{T}\frac{\tilde{\psi}^{''}(Y_{s })}{2}(|m_{s}Z_{s}|^{2}dC_{s} + d\lan{L} \ran_{s}) \right) .} $ \\
\end{tabbing}
Since $Z \cdot M $ and $L$ are square integrable martingales and $\tilde{\psi}^{'}(Y)$ is a bounded process, the conditional expectation of the terms of the second line in the right-hand side vanishes. 
Using both the upper bound on $F$ in ($H_{1}$) and simple computations, we obtain\\
\begin{tabbing}
 $ $ \= $\disp{\tilde{\psi}(Y_{\tau}) - \mathbb{E}^{\mathcal{F}_{\tau }}(\tilde{\psi}(Y_{T}))} $ \\
$\quad \quad$ \> $\disp{ \leq \mathbb{E}^{\mathcal{F}_{\tau}}\int_{\tau}^{T}\tilde{\psi}^{'}(Y_{s })(|\bar{\alpha}_{s}|(1 + b|Y|_{S^{\infty}})dC_{s}} $ \\
 
\> $ \quad \quad \disp{+ \;\mathbb{E}^{\mathcal{F}_{\tau }}\int_{\tau}^{T}(\frac{\beta}{2}\tilde{\psi}^{'} - \frac{1}{2} \tilde{\psi}^{''})(Y_{s})d\lan{L} \ran_{s}}
 \disp{+ \;\mathbb{E}^{\mathcal{F}_{\tau}}\int_{\tau}^{T}(\frac{\gamma}{2}\tilde{\psi}^{'} - \frac{1}{2} \tilde{\psi}^{''})(Y_{s})|m_{s}Z_{s}|^{2}dC_{s}}. $ \\
\end{tabbing}
\indent Putting in the left-hand side of this formula the two last terms and 
using the properties of $\psi_{\gamma}$ given by (\ref{eq: Proppsi}) and the fact that: $\gamma \geq |\beta|$,
we get
$$ \big(\frac{1}{2} \tilde{\psi}^{''}- \frac{\beta}{2}\tilde{\psi}^{'}\big)(Y_{s}) \geq \big(\frac{1}{2} \tilde{\psi}^{''}- \frac{\gamma}{2}\tilde{\psi}^{'}\big)(Y_{s}) =  \frac{1}{2},\; \; \mb{P}\textrm{-a.s. and for all}\; s.  $$
It then follows from these two last inequalities 
 \begin{tabbing}
$\disp{ \mathbb{E}^{\mathcal{F}_{\tau } }  \left(\int_{\tau }^{T}|m_{s}Z_{s}|^{2}dC_{s} + (\lan{L} \ran_{T } - \lan{L} \ran_{\tau}) \right)}$\\
\\
 $\quad \leq 
\disp{2\mathbb{E}^{\mathcal{F}_{\tau } }\left( \int_{\tau}^{T}(\frac{1}{2} \tilde{\psi}^{''}- \frac{\gamma}{2}\tilde{\psi}^{'}\big)(Y_{s})|m_{s}Z_{s}|^{2} dC_{s}
\; + \int_{\tau}^{T}(\frac{1}{2} \tilde{\psi}^{''}- \frac{\beta}{2}\tilde{\psi}^{'}\big)(Y_{s}) d\lan{ L \ran}_{s}\right) \leq C^{'}},$\\
 \end{tabbing} 
with a constant $C^{'}$ independent of the stopping time $\tau$ and depending only on the parameters $a$, $b$, $\gamma$ and $|B|_{\infty}$.
\begin{flushright}
 $\square$
\end{flushright}


\subsection{The uniqueness result}
\textbf{Proof:}
The key idea of this proof is to proceed by linearization and to justify as in \cite{ImkelleretHu} the use of Girsanov's theorem.
 Let ($Y^{1}$, $Z^{1}$, $L^{1}$) and ($Y^{2}$, $Z^{2}$, $L^{2}$) be two solutions of the BSDE($F, \beta, B$) with $F$ satisfying both $(H_{1})$ and $(H_{2})$ and $B$ bounded. 
 We define $Y^{1, 2}$ by: $Y^{1, 2} =Y^{1} - Y^{2}$ ($Z^{1, 2}$ and $L^{1, 2}$ are defined similarly)
and we consider the nonnegative and bounded
semimartingale ($\tilde{Y}^{1, 2}$) such that, for all $t$,  
$ \tilde{Y}^{1, 2}_{t} = e^{2\mu C_{t}}|Y_{t}^{1, 2}|^{2} $.
We then use It\^o's formula  
$$d\tilde{Y}_{s}^{1, 2}  = 2\mu \tilde{Y}_{s}^{1, 2}dC_{s}  
 + e^{2\mu C_{s}}2Y_{s}^{1, 2}dY_{s}^{1, 2}
+ \frac{1}{2}e^{2\mu C_{s}}2 d\lan{Y}^{1, 2}\ran_{s}.$$
$Y^{1}$ and $Y^{2}$ being solutions of the BSDE($F, \beta, B$),   $$ dY_{s}^{1, 2} = -(F(s, Y_{s}^{1}, Z_{s}^{1}) - F(s, Y_{s}^{2}, Z_{s}^{2}))dC_{s} -\frac{\beta}{2}d(\lan{L}^{1}\ran_{s} - \lan{L}^{2}\ran_{s}) + dK_{s}, $$  with: $K = \; Z^{1, 2} \cdot M + L^{1, 2}$, which stands for the martingale part. Hence,
integrating It\^o's formula between $t$ and an arbitrary $\mc{F}$-stopping time $\tau$\\
 \[ \begin{array}{ll}
 \disp{  \tilde{Y}_{t}^{1, 2}  - \tilde{Y}_{\tau}^{1, 2} = }  &
 \disp{ - \int_{t}^{\tau} 2\mu \tilde{Y}_{s}^{1, 2}dC_{s} }\\
 & \disp{  + \int_{t}^{\tau}e^{2\mu C_{s}}2Y_{s}^{1, 2}(F(s, Y_{s}^{1}, Z_{s}^{1}) - F(s, Y_{s}^{2}, Z_{s}^{2}))dC_{s}}  \\ 
& \disp{ + \int_{t}^{\tau}e^{2\mu C_{s}}2Y_{s}^{1, 2}\frac{\beta}{2}d\lan{L}^{1,2}, L^{1}+ L^{2}\ran_{s} }  \\
& \disp{ - \int_{t}^{\tau}e^{2\mu C_{s}}2Y_{s}^{1, 2}\big(Z_{s}^{1, 2}dM_{s} + d L_{s}^{1, 2} \big)}  \\
&\disp{  \underbrace{-\int_{t}^{\tau }e^{2\mu C_{s}}\frac{1}{2}2d\lan{Y}^{1, 2}\ran_{s}}_{ \leq 0}} \; .\\
\end{array} \] 
 $F$ satisfying ($H_{2}$), it follows
\[  2Y_{s}^{1, 2}(F(s, Y_{s}^{1}, Z_{s}^{1}) - F(s, Y_{s}^{2}, Z_{s}^{2})) \leq  2\mu |Y_{s}^{1, 2}|^{2} + 2Y_{s}^{1, 2}(m_{s}\kappa_{s})^{'}(m_{s}Z_{s}^{1, 2}) , \] 
where the $\mb{R}^{d}$-valued process $\kappa$ is defined as 
\[ \left\{ \begin{array}{ll}
 \disp{\kappa_{s} = \frac{(F(s, Y_{s}^{2}, Z_{s}^{1}) - F(s, Y_{s}^{2}, Z_{s}^{2}))(Z_{s}^{1, 2})}{|m_{s}(Z_{s}^{1, 2})|^{2}}}, & \textrm{if } \; |m_{s}(Z_{s}^{1, 2})| \neq 0, \\
 \disp{\kappa_{s} = 0,} & \textrm{otherwise.} \\
\end{array} \right.
 \]
We introduce a new process $A$
\[ \begin{array}{ll}
 A_{s} := & \big( 2Y_{s}^{1, 2}(F^{1}(s, Y_{s}^{1}, Z_{s}^{1}) - F^{2}(s, Y_{s}^{2}, Z_{s}^{2})\big)\\& - \left(2\mu|Y_{s}^{1, 2}|^{2} + 2Y_{s}^{1, 2}(m_{s}\kappa_{s})^{'}(m_{s}Z_{s}^{1, 2}) \right).\\
\end{array} \]
The process $A$ being almost surely non positive, we obtain\\
\[ \begin{array}{l}
\disp{ \tilde{Y}_{t}^{1, 2} -\tilde{Y}_{\tau}^{1, 2}  }  
\quad \; \disp{= \underbrace{\int_{t}^{\tau}A_{s}dC_{s}  -\int_{t}^{\tau}e^{2\mu C_{s}}\frac{1}{2}2d\lan{Y}^{1, 2}\ran_{s}}_{\le 0}}\\

 \; \disp{\; + \int_{t}^{\tau}2Y_{s}^{1, 2}e^{2 \mu C_{s}}(m_{s}\kappa_{s})^{'}(m_{s}Z_{s}^{1, 2})dC_{s}} \\
\disp{ \;+\; \int_{t}^{\tau}2Y_{s}^{1, 2}e^{2 \mu C_{s}}\frac{\beta}{2}d\lan{L}^{1, 2}, L^{1} + L^{2}\ran_{s}}  \\
 \; \disp{ \; -
   \int_{t}^{\tau}2e^{2 \mu C_{s}}Y_{s}^{1, 2}Z_{s}^{1, 2}dM_{s}  -  \int_{t}^{\tau}2e^{2 \mu C_{s}}Y_{s}^{1, 2}dL_{s}^{1, 2} }. \\
 \end{array} \]
We then consider the following stochastic integrals  
\[ \begin{array}{ll}  \disp{\tilde{N} = \left(2e^{2\mu C_{t}}Y^{1,
  2}Z^{1, 2}\right) \cdot M \; \textrm{and} \;  \bar{N} = \kappa \cdot M, } & \textrm{ on the one hand,}\\
 \disp{\tilde{L} =  \left(2Y^{1,
  2}e^{2\mu C}\right) \cdot L^{1, 2} \; \textrm{and} \; \bar{L} = \frac{\beta}{2}(L^{1} + L^{2}),} &\; \textrm{ on the other hand.} \\
\end{array} \]
 From ($H_{2}$), we deduce  
$$|m_{s}\kappa_{s}| \leq C (|m_{s}\theta_{s}| + |m_{s}Z_{s}^{1}| + |m_{s}Z_{s}^{2}|).$$
Using both the assertion (ii) in Lemma \ref{estimapriori} and the assumption on $\theta$ given by ($H_{2}$), we get:  $\kappa \cdot M + \frac{\beta}{2}(L^{1} + L^{2})$ is a BMO martingale. Hence, by referring to \cite{Kazamaki}, $\mathcal{E}(\kappa \cdot M + \frac{\beta}{2}(L^{1} + L^{2}))$
is a true martingale. Defining $\mathbb{Q}$ such that: $d\mathbb{Q} = \mathcal{E}(\kappa \cdot M + \frac{\beta}{2}(L^{1} + L^{2}))d\mathbb{P}$, Girsanov's theorem entails that $K :=\tilde{N} + \tilde{L} - \lan{\tilde{N}} + \tilde{L}, \kappa \cdot M + \frac{\beta}{2}(L^{1} + L^{2})\ran$ is a local martingale under $\mathbb{Q}$. Hence, taking the conditional expectation w.r.t $\mc{F}_{t} $ and under $\mb{Q} $, we get
  
$$\disp{ \tilde{Y}_{t}^{1, 2}  \le  \mb{E}^{\mb{Q}}\left(\tilde{Y}_{\tau}^{1, 2} | \mc{F}_{t}\right) } . $$
 The same localization procedure as in the proof of lemma \ref{estimapriori} at the top of page 10 (replacing the word submartingale by supermartingale and using then the bounded convergence theorem) leads to
$$\forall \; t, \quad \tilde{Y}_{t}^{1,2} \leq 0 \; \;
  \mathbb{Q}\textrm{-a.s.}
\quad \textrm{(and } \; \mathbb{P}\textrm{-a.s., because of the equivalence of}\; \mathbb{P} \; \textrm{and} \; \mathbb{Q}), $$
which ends the proof ($\tilde{Y}^{1,2} $ being a non negative process).\\

\subsection{Existence }

\subsubsection{Main steps of the proof of Theorem \ref{existence}}

\hspace{0.5cm} In this part and to establish the existence result (Theorem \ref{existence}), we proceed with three main steps.\\
 \indent In a first step, we prove that, to solve a BSDE of type ($\textrm{Eq1}$) under assumption ($H_{1}$), it suffices to solve the same BSDE under a simpler assumption ($H_{1}^{'}$).\\
\indent In a second step, we introduce an intermediate BSDE of the form ($\textrm{Eq2}$) and we establish a one to one correspondence between the existence of a solution of a BSDE of the form ($\textrm{Eq1}$) and one of the form ($\textrm{Eq2}$).\\
\indent The third and last step consists in constructing a solution of the BSDE of the form ($\textrm{Eq2}$) when its generator $g$ satisfies ($H_{1}^{'}$) using a ``monotone stability'' result analogous to the one given in \cite{mkobylanski}.\\ 

\paragraph{Step 1: Truncation in $\mathbf{y}$}

\indent We rely on the a priori estimates given in Lemma \ref{estimapriori} to strengthen the assumption on the generator and obtain precise estimates for an intermediate BSDE.
More precisely, we restrict ourselves to the simpler assumption ($H_{1}^{'}$) (instead of ($H_{1}$))
\begin{eqnarray*}
 (H_{1}^{'}) \quad \exists \; \bar{\alpha} \; \geq 0 \; \int_{0}^{T}\bar{\alpha}_{s}dC_{s} \leq a \; (a > 0), \; \textrm{s.t.} \;  |F(s, y, z)| \leq \bar{\alpha}_{s} + \frac{\gamma}{2}|m_{s}z|^{2} .
\end{eqnarray*}
Assuming here that we have a solution of the BSDE($F, \beta, B$) of type ($\textrm{Eq1}$) under assumption ($H_{1}^{'}$) on $F$, we deduce the existence of a solution of this BSDE under ($H_{1}$). For this, we define $K$ by: $K$ = $|c| + |C|$, with $c$ and $C$ the two constants given in (i), Lemma 1 and we introduce 
\[ \left\{ \begin{array}{l}
dY_{s}^{K} = -F^{K}(s, Y_{s}^{K}, Z_{s}^{K} )dC_{s} -\frac{\beta}{2}d\lan{L}^{K}\ran_{s} + Z_{s}^{K} dM_{s} + dL_{s}^{K}, \\
Y_{T}^{K} = B,
 \end{array} \right. \]
where $F^{K}$ and the truncation function $\rho_{K}$ are respectively defined by:
 $F^{K}(s, y, z) = F(s, \rho_{K}(y), z)$ and
\[ \rho_{K}(x) =  \left\{  \begin{array}{ll} 
 - K &\textrm{if} \; x < - K ,\\
 x &\textrm{if} \; |x| \leq K ,\\
 K &\textrm{if} \; x > K . \\
 \end{array} \right. \]
Hence, we have
$$ \forall \; y \in \mathbb{R}, \; z \in \mathbb{R}^{d}, \; |F^{K}(s, y, z)|  \leq \bar{\alpha}_{s}(1 + b|\rho_{K}(y)|) + \frac{\gamma}{2}|m_{s}z|^{2} .$$
Since: $|\rho_{K}(x)| \leq |x|$, $F^{K}$ satisfies again ($H_{1}$) with the same parameters as $F$. Using Lemma \ref{estimapriori}, $K$ is an upper bound of $Y^{K}$ in $S^{\infty}$, for any solution ($Y^{K}, Z^{K}, L^{K})$ of BSDE($F^{K}$, $\beta$, $B$). Besides, if we replace $\bar{\alpha}$ by $\tilde{\alpha} := \bar{\alpha}(1 + b K)$, $F^{K}$ satisfies ($H_{1}^{'}$).
Due to the initial assumption, there exists a solution denoted by ($Y^{K}$, $Z^{K}$, $L^{K}$) of BSDE($F^{K}$, $\beta$, $B$). Since: 
$|Y^{K}| \leq K$, $F^{K}$ and $F$ coincide along the trajectories of this solution and hence, ($Y^{K}$, $Z^{K}$, $L^{K}$) is a solution of BSDE($F$, $\beta$, $B$) with $F$ satisfying $(H_{1})$.\\

\paragraph{Step 2: an intermediate BSDE} 
\indent To establish the one-to-one correspondence,
we first assume the existence of a solution ($Y$, $Z$, $L$) of BSDE($F, \beta, B$) with $F$ satisfying ($H_{1}^{'}$) and we set: $U:=e^{\beta Y}$.
Thanks to It\^o's formula, $U$ solves a
BSDE of the type ($\textrm{Eq2}$) and the expression of the generator $g$ is 
\begin{equation}\label{eq: generateurg}
g(s, u, v) = \big(\beta u F(s, \frac{\ln(u)}{ \beta},  \frac{v}{
  \beta u}) -  \frac{1}{2u}|m_{s}v|^{2}\big)\mathbf{1}_{u >0},  \end{equation}
and a solution of the BSDE($g, e^{\beta B}$) of type ($\textrm{Eq2}$) is a triple ($U$, $V$, $N$) such that: $U_{s}:= e^{\beta Y_{s}}$,
$V_{s} = \beta U_{s}Z_{s}$ and $N = \beta U \cdot L$.
Our aim is to prove that the converse is true: i.e. if we can solve the BSDE($g$, $e^{\beta B}$) of type ($\textrm{Eq2}$) under the assumption ($H_{1}^{'}$) on $g$, then we obtain a solution of the BSDE($F$, $\beta$, $B$) of type ($\textrm{Eq1}$) by setting\\
\begin{eqnarray}\label{eq: formulcorrespondance}
Y: =  \frac{\ln(U)}{ \beta}, \; Z: = \; \frac{V}{\beta U},\;  \textrm{and}\; L: =
\frac{1}{\beta U} \cdot N.
\end{eqnarray}
\indent To achieve this, we give precise estimates of $U$ in $S^{\infty}$ for any solution ($U$, $V$, $N$) of the BSDE($g$, $e^{\beta B}$) of type ($\textrm{Eq2}$). Due to the singularity of the expression (\ref{eq: generateurg}) of $g$ with respect to $u$, we first rely on a truncation argument and therefore, we introduce a new generator $G$ 
 \[ G(s, u ,v) = \beta \rho_{c^{2}}(u)F(s, \frac{\ln(u \vee c^{1})}{ \beta},  \frac{v}{
  \beta(u \vee c^{1})})-  \frac{1}{2(u \vee c^{1})}|m_{s}v|^{2}.  \]
The two positive constants $c^{1}$ and $c^{2}$ are defined later and the function $\rho_{c^{2}}$ is the same as in the first step. $F$ satisfying ($H_{1}^{'}$) and since: $\rho_{c^{2}}(u) \leq c^{2}$, we obtain that $G$ also satisfies ($H_{1}^{'}$). Hence, for any $c^{1}$, $c^{2}$, there exists a solution of the BSDE($G, e^{\beta B}$) of type ($\textrm{Eq2}$). We denote it ($U^{c^{1},\;c^{2}}$, $V^{c^{1},\;c^{2}}$, $N^{c^{1},\;c^{2}}$).
Thanks to the estimates
\[ \begin{array}{ll} |G(s, u ,v)| &\leq \beta \rho_{c^{2}}(u)( \bar{\alpha_{s}} + \frac{\gamma|m_{s}v|^{2}}{2|\beta c^{1}|^{2}}) + \frac{|m_{s}v|^{2}}{2c^{1}}\\

 &\leq \beta \bar{\alpha_{s}}|u| + \frac{\hat{\gamma}}{2}|m_{s}v|^{2}, \; \textrm{with} \; \hat{\gamma} := \frac{\gamma c^{2}}{|\beta||c^{1}|^{2}} + \frac{1}{c^{1}},\\
    \end{array}\]
 $G$ satisfies ($H_{1}$) with parameters $a, b, \gamma$ such that $$ a: =  \;\disp{\int_{0}^{T}|\beta|\bar{\alpha}_{s}dC_{s}},\;  b: = 1,  \; \gamma:= \hat{\gamma}.$$ Using (i) in Lemma \ref{estimapriori}, the solution ($U^{c^{1},\;c^{2}}$, $V^{c^{1},\;c^{2}}$, $N^{c^{1},\;c^{2}}$) satisfies \\ $$ U^{c^{1},\;c^{2}} \leq e^{a} - 1 + |e^{\beta B}|_{\infty}e^{a},\;  \mathbb{P}\; \textrm{-a.s.} $$ 
 Defining $c^{2}$ by: $c^{2}$ := $e^{a} - 1 + |e^{\beta B}|_{\infty}e^{a}$, this provides an upper bound independent of $\gamma$. To prove the existence of a strictly positive lower bound, we consider a solution ($U$, $V$, $N$) of the BSDE($G$, $e^{\beta B}$) and
 we define the adapted process $\Psi(U)$ for all $t$ by 
$  \Psi(U_{t}) = e^{-\int_{0}^{t}\tilde{\beta}_{s}dC_{s}}U_{t} $ (
$\tilde{\beta} := |\beta|\bar{\alpha}\textrm{sign}(U_{s})$) is such that: 
 $\disp{\int_{0}^{T}|\tilde{\beta}_{s}|dC_{s}} \leq a$, $\mb{P}$-a.s.).
Applying then It\^o's formula to $ \Psi(U)$ between $t$ and $T$, we get 
\begin{tabbing}
$ $ \= $ \disp{\Psi(U_{t}) -} \disp{ \Psi(U_{T}) }$ \\
\\
 \> $ \;= \; \disp{ \int_{t}^{T}\big( e^{-\int_{0}^{s}\tilde{\beta}_{u}dC_{u}}( G(s, U_{s}, V_{s}) + \tilde{\beta}_{s}U_{s})\big)dC_{s}  - \int_{t}^{T} e^{-\int_{0}^{s}\tilde{\beta}_{u}dC_{u}}(V_{s}dM_{s} + dN_{s}) }$.\\
\> $\;=\; \disp{ \int_{t}^{T}e^{- \int_{0}^{s}\tilde{\beta}_{u}dC_{u}}A_{s}dC_{s} - \int_{t}^{T}\frac{\gamma}{2}e^{- \int_{0}^{s}\tilde{\beta}_{u}dC_{u}}|m_{s}V_{s}|^{2}dC_{s} }$\\
\> $ \quad \; \; 
 - \;\disp{\int_{t}^{T}\big(e^{- \int_{0}^{s}\tilde{\beta}_{u}dC_{u}}(V_{s}dM_{s} + dN_{s})\big)}, $
\end{tabbing}
with the process $A$ such that: $A_{s} := G(s, U_{s}, V_{s}) + \big(\tilde{\beta}_{s}U_{s} + \frac{\gamma}{2}|m_{s}V_{s}|^{2}\big)$ which is almost surely positive. 
Since $ -\frac{\gamma}{2}(V \cdot M)$ is a BMO martingale (thanks to (ii) in Lemma \ref{estimapriori}), we introduce a probability measure by defining: $\frac{d\mathbb{Q}}{d\mathbb{P}} = \mathcal{E}(-\frac{\gamma}{2}V \cdot M)$. The Girsanov's transform $\tilde{M}$ of $M$ such that:
$\tilde{M}  := \; M + \frac{\gamma}{2}\lan{V \cdot M},M\ran$, is a local martingale under $\mathbb{Q}$ and it follows that $\Psi(U)$ is a local submartingale under $\mb{Q}$: relying now on the standard localization procedure and on the boundeness assumption on $\Psi(U)$, we conclude
$$  \Psi(U_{t}) \geq  \mathbb{E}^{\mathbb{Q}}(\Psi(U_{T})| \mathcal{F}_{t}). $$
Hence, $U_{t} \geq \mathbb{E}^{\mathbb{Q}} (\big(\mathop{\inf{U_{T}}}\big)e^{-\int_{t}^{T}\tilde{\beta}_{s}dC_{s}} | \mathcal{F}_{t})$, and
if $c^{1}$ is defined by $c^{1}: = e^{-|\beta| \big( |B|_{\infty} + a \big)}$, it is a lower bound of $U$.
For these choices of $c^{1}$, $c^{2}$, 
 the generator $G$ satisfies ($H_{1}$) and, for any solution ($U$, $V$, $N$),
 $$    c^{1} \leq U_{s} \leq c^{2}, \; \;\mathbb{P}\textrm{-a.s. and for all} \; \textrm{s}. $$
Since: $G(s, U_{s}, V_{s}) = g(s, U_{s}, V_{s}) \; \mathbb{P}\textrm{-a.s.}$ and for all $s$, ($U$, $V$, $N$) is a solution of the BSDE($g$, $e^{\beta B}$).
The process $U$ being strictly positive and bounded, we can define ($Y$, $Z$, $L$) by using (\ref{eq: formulcorrespondance}) and, applying It\^o's formula to $\frac{\ln(U)}{\beta}$, we check that ($Y$, $Z$, $L$) is a solution of the BSDE($F, \beta, B$).\\

\paragraph{Step 3: Approximation}
\indent To prove the existence of a solution of the BSDE($F, \beta, B$) of type ($\textrm{Eq1}$) under ($H_{1}$),
the above two steps show that it is sufficient to prove the \\ existence of a solution of the BSDE($g, e^{\beta B}$) of type ($\textrm{Eq2}$) under ($H_{1}^{'}$) on $g$.
Analogously to \cite{mkobylanski}, we construct an approximating sequence ($U^{n}$, $V^{n}$, $N^{n}$) satisfying\\
$\bullet \; $ these triples are solutions
of the BSDEs($g^{n}$, $e^{\beta B}$),\\
$\bullet \;$ the sequence ($g^{n}$) is increasing and converges, $\mb{P}$-a.s. and for all $s$, to $g \;(g: (y, z) \to g(s, y, z))$.\\
\indent From now and for the remaining of Section 3.3.1, we suppose\\
\begin{equation}\label{eq: assumpt1}
\textrm{\textit{Assumption 1:} The driver $g$ satisfies ($H_{1}^{''}$).}
\end{equation}
 We then proceed by defining $g^{n}$ by inf-convolution
\[\; g^{n}(s, u, v) =  \textrm{ess}\displaystyle{\inf_{u^{'},\; v^{'}, \atop\; u^{'}, v^{'} \in \mathbb{Q}^{d} }\left(g(s, u^{'},v^{'} )  + n|m_{s}(v -  v^{'})|\; +\; n|u -u^{'}| \right) . }   \] 
Such a $g^{n}$ is well defined and globally Lipschitz continuous in the following sense
\begin{equation}\label{eq: conditionlip}
\forall \; u^{1}, u^{2}, \; v^{1}, v^{2}, \quad |g^{n}(s, u^{1}, v^{1}) - g^{n}(s, u^{2}, v^{2})| \leq n\big(|m_{s}(v^{1} - v^{2})| +\; |u^{1} -u^{2}|\big). \end{equation}
Since ($g^{n}$) is increasing and converges, $\mb{P}$-a.s. and for all $s$, to $g:(u,v) \to g(s,u,v)$ which is continuous w.r.t. ($u, v$), Dini's theorem implies that the convergence is uniform over compact sets.
Besides using that: $g^{n} \le g$, we obtain 
\begin{equation}\label{eq: condition2}  \displaystyle{\sup_{n}|g^{n}(s, 0, 0)|} \leq \bar{\alpha}_{s}.
\end{equation}
\indent The existence of a unique solution ($U^{n}$, $V^{n}$, $N^{n}$) of the BSDEs given by ($g^{n}$, $e^{\beta B}$) in $S^{2} \times L^{2}(d\lan{M} \ran \otimes d\mathbb{P}) \times \mathcal{M}^{2}([0,T])$ follows from (\ref{eq: conditionlip}) and (\ref{eq: condition2}) (we refer to \cite{ElkaretMazliak} for results in a general filtration). $S^{2}$ denotes here the space of all continuous processes $U$ such that: $\mb{E}\big(\disp{\sup_{t \in [0, T]}|U_{t}|^{2}} \big) < \infty$. Furthermore, applying Theorem \ref{comparison} for these BSDEs of type ($\textrm{Eq2}$) and using that $(g^{n})_{n}$ is increasing, we get: $U^{n} \le U^{n+1}$. The following result entails that, for all $n$, $U^{n}$ is in $S^{\infty}$.

\begin{proposition}\label{estimclass}
Let ($U^{n}$, $V^{n}$, $N^{n}$) be a solution in $S^{2} \times L^{2}(d\lan{M} \ran \otimes d\mathbb{P}) \times \mathcal{M}^{2}([0,T])$ of a BSDE of the type ($\textrm{Eq2}$) given by the parameters ($g^{n}$, $\bar{B}$), with a generator $g^{n}$ $L_{n}$-Lipschitz and a terminal condition $\bar{B}$ bounded, we have 
\begin{equation}\label{eq: controlbornee}
\exists \; K(L_{n}, T) > 0, \;\forall \;t, \; |U_{t}^{n}|^{2} \leq K(L_{n}, T)\mb{E}\left(|\bar{B}|^{2} + (\int_{t}^{T}|g^{n}(s, 0, 0)|dC_{s})^{2}| \mc{F}_{t} \right). \end{equation}
\end{proposition}
The proof, relegated to the appendix, is adapted from the results given in Proposition 2.1 in \cite{BriandetCoquet}. Relying on (\ref{eq: condition2}) and on the assumption on $\bar{\alpha}$, Proposition \ref{estimclass} implies that $U^{n}$ is in $S^{\infty}$. Now, since each generator $g^{n}$ satisfies the assumption ($H_{1}^{''}$) (and hence ($H_{1}$) with the same parameters), we can use assertion (i) in Lemma \ref{estimapriori} to ensure that ($U^{n}$) is uniformly bounded in $S^{\infty}$. \\
\newline
\paragraph{Step 4: Convergence of the approximation}
 To prove the convergence of the solutions of the BSDEs($g^{n}, e^{\beta B}$) under \textit{Assumption 1} (see (\ref{eq: assumpt1})), we introduce the triplet ($\tilde{U}, \tilde{V}, \tilde{N}$) as being the limit (in a specific sense) of ($U^{n}, V^{n}, N^{n}$).
($U^{n}$) being increasing, we set:
$\tilde{U}_{s} = \displaystyle{\lim_{n} \nearrow (U_{s}^{n})}, \; \mb{P}\textrm{-a.s.} \;$ and for all $s$.
Any generator $g^{n}$ satisfying ($H_{1}^{''}$), and hence ($H_{1}$) with the same parameters, the estimate (ii) in Lemma \ref{estimapriori} holds true for each term of ($V^{n}$)$_{n}$ and ($N^{n}$)$_{n}$ (uniformly in $n$).
As bounded sequences of Hilbert spaces, there exist subsequences of $(V^{n})$ and $(N_{T}^{n})$ such that
 $V^{n} \xrightarrow{w} \tilde{V}$ (in $L^{2}(d\lan{M} \ran \otimes d\mathbb{P})$), and $N_{T}^{n}  \xrightarrow{w} \tilde{N}_{T}$ in $L^{2}(\Omega, \mc{F}_{T}, \mb{P})$. This implies the weak convergence in $L^{2}(\Omega, \mc{F}_{t}, \mb{P})$ of $N_{t}^{n}$ to $\tilde{N}_{t}$, if we define $\tilde{N}_{t}$ by: $\tilde{N}_{t} := \mathbb{E}^{\mathcal{F}_{t}}(\tilde{N}_{T})$.
However, to justify the passage to the limit in the BSDEs given by ($g^{n}$, $e^{ \beta B}$), we need the strong convergence of ($V^{n}$), eventually along a subsequence, to $\tilde{V}$ in $L^{2}(d\lan{M} \ran \otimes d\mathbb{P}) $ (resp. ($N^{n}$) to $\tilde{N}$ in $ \mc{M}^{2}([0, T])$).
We give one essential result (similar to the stability result in \cite{mkobylanski}) which is the key ingredient in the last step of the proof of Theorem \ref{existence}.

 \begin{proposition}\label{convmonotone}
Let ($g^{n}$) and ($\tilde{B}^{n}$) be two sequences associated with the BSDEs($g^{n}, \tilde{B}^{n}$) of type ($\textrm{Eq2}$) and satisfying \\
 $\bullet \;$ $\mb{P}$-a.s. and for all $s$, ($g^{n}:(u, v) \to g^{n}(s, u, v)$) converges increasingly w.r.t. $n$ and uniformly on the compact sets of $\mathbb{R} \times \mathbb{R}^{d}$ to
$g$ ($g:(u, v) \to  g(s, u, v)$ ($g$ is continuous w.r.t. ($u, v$)).\\
  $\bullet \;$ For all $n$, each $g^{n}$ satisfies ($H_{1}^{''}$), with the same parameters as $g$ (independent of $n$),\\
$\bullet \;$ $(\tilde{B}^{n})$ is a uniformly bounded sequence of $\mathcal{F}_{T}$-measurable random variables, which converges almost surely to $\tilde{B}$ and increasingly w.r.t. $n$.\\
If there exists one solution ($U^{n}$, $V^{n}$, $N^{n}$) of the BSDEs given by ($g^{n}$, $\tilde{B}^{n}$) such that the sequence ($U^{n}$)$_{n}$ is increasing, then the sequence ($U^{n}$, $V^{n}$, $N^{n}$) converges to ($\tilde{U}, \tilde{V}, \tilde{N}$) in the following sense
\[ \textrm{and} \left\{ \;\begin{array}{l}
 \disp{\mb{E}(\sup_{t \in [0, T]} |U_{t}^{n} - \tilde{U}_{t}|)} \to 0, \;\textrm{as} \; n \; \to \infty, \\
   \mathbb{E}\left(\disp{\int_{0}^{T}|m_{s}(\tilde{V}_{s} -
   V_{s}^{n})|^{2}dC_{s} + |\tilde{N}_{T} - N_{T}^{n}|^{2}}  \right) \to 0, \;\textrm{as} \; n \; \to \infty. \\
   \end{array} \right. \]
Besides, ($\tilde{U}, \tilde{V}, \tilde{N}$) is solution of the BSDE($g$, $\tilde{B}$) of type ($\textrm{Eq2}$).
 \end{proposition}
\textbf{Remark:}
This ``stability'' result stated in Lemma \ref{convmonotone} holds also for the solution of the BSDE($F, \beta, B$) of type ($\textrm{Eq1}$) and the proof is obtained using the correspondence established in the second step.\\

 We relegate to subsection 3.3.2 the technical point in the proof of Lemma \ref{convmonotone}, i.e. the strong convergence in their respective Hilbert spaces of the sequences ($V^{n}$) and ($N^{n}$). Assuming this result, we prove the existence of a solution for BSDE($g, \; \tilde{B}$) by justifying the
passage to the limit in BSDEs($g^{n}, \tilde{B}^{n}$)\\
\[  U_{t}^{n} = \tilde{B}^{n} + \int_{t}^{T} g^{n}(s, U_{s}^{n}, V_{s}^{n})dC_{s} - \int_{t}^{T}V_{s}^{n}dM_{s} -(N_{T}^{n} - N_{t}^{n}). \]
To this end, we check that, $\mathbb{P}$-a.s. and for all $t$,  \\
 (i) $V^{n} \to
 \tilde{V} \; \; (\textrm{in} \; L^{2}(d\lan{ M\ran} \otimes d\mb{P}))$, as $n \to \infty$,\\
  (ii) $\disp{N^{n} \rightarrow \tilde{N} \; (\textrm{in} \;\mc{M}^{2}([0, T]))} $, as $n \to \infty$,\\
  (iii)  $\disp{\mb{E}\big(\int_{0}^{t}|g^{n}(s,U_{s}^{n},V_{s}^{n}) - g(s,\tilde{U}_{s},\tilde{V}_{s})|dC_{s}\big) \to 0 }$, as $n \to \infty$.\\
Assertions (i) and (ii) are consequences of 
the strong convergence of the
sequences ($V^{n}$) (resp. ($N^{n}$)) in $L^{2}(d \lan{M}\ran \times d\mb{P})$ (resp. in $\mc{M}^{2}([0, T])$).
To prove (iii), we justify the convergence in $L^{1}(ds \otimes d\mb{P})$ using the two following results:\\
 $\bullet \;$ The convergence in $dC_{s} \otimes d\mb{P}$-measure of ($m_{s} V_{s}^{n}$) and ($U_{s}^{n}$) (at least along proper subsequences) and the properties of ($g^{n}$), which ensure the convergence of ($g^{n}(s,U_{s}^{n},V_{s}^{n}))$ to $g(s,\tilde{U}_{s},\tilde{V}_{s}) $ in $dC_{s} \otimes d\mb{P}$-measure.\\
$\bullet \;$  The uniform integrability of the family ($g^{n}(s,U_{s}^{n},V_{s}^{n})$) resulting from the estimates of $g^{n}$ given by ($H_{1}^{'}$) and from the fact that ($|m V^{n}|^{2}$) is a uniformly integrable sequence, since it is strongly convergent in $L^{1}(dC \times d\mb{P})$. \\
Passing to the limit as $n$ goes to $\infty$, we get that the triplet ($\tilde{U}, \tilde{V}, \tilde{N}$)  
is a solution of the BSDE($g$, $e^{ \beta B}$).\\
To obtain a solution of the BSDE($F, \beta, B$), we rely on the results of the two first steps and we set ($\tilde{Y}$, $\tilde{Z}$, $\tilde{L}$) using the formula (\ref{eq: formulcorrespondance}).
\\
\begin{flushright}
$   \square $
\end{flushright}
\indent 
Now, we relax \textit{Assumption 1} given by (\ref{eq: assumpt1}): i.e., we proceed with the case when $g$ only satisfies $(H_{1}^{'})$. In this case, the lower bound is no more Lipschitz and, for the method, we refer once again to \cite{BriandetHu}: the idea consists in using two successive approximations. For this, we define ($g^{n,\; p}$) as follows
\begin{tabbing}  
$g^{n, p}(s, u, v) =  $\= $ \; \textrm{ess} \displaystyle{\inf_{ u^{'}, v^{'}}\big(g^{+}(s, u^{'}, v^{'}) + n|m_{s}(v-v^{'})| +\; n|u -u^{'}|\big)} $\\  \> $\quad  - \; \textrm{ess} \displaystyle{\inf_{u^{'}, v^{'}}\big(g^{-}(s, u^{'}, v^{'}) + p|m_{s}(v-v^{'})| +\; p|u -u^{'}|\big) ,}$ \\ \end{tabbing}
which is increasing w.r.t. $n$ and decreasing w.r.t. $p$. The entire proof can be rewritten identically by passing to the limit as $n$ goes to $\infty$ ($p$ being fixed) and then as $p$ goes to $\infty$.  

\subsubsection{Proof of the ``stability'' result in Lemma \ref{convmonotone}}
\hspace{0.5cm} Following the same method as in \cite{mkobylanski}, we establish the strong convergence of the sequences
($V^{n} $)$_{n}$ and ($N^{n}$)$_{n}$ to $\tilde{V}$ and
$\tilde{N}$ (this requires the a priori estimates established in
Lemma 1 for the solutions of the BSDEs given by ($g^{n}$, $\tilde{B}^{n} $)). We first introduce the nonnegative semimartingale $\Phi_{L}(U^{n} -U^{p})$ = $(\Phi_{L}(U^{n, p}))_{ n \geq p}$, 
with $\Phi_{L}$ such that \\
\begin{equation}\label{eq: PhiL} \Phi_{L}(x) = \frac{e^{Lx}- Lx - 1 }{L^{2}}. \end{equation}
$\Phi_{L}$ satisfies: $\Phi_{L}\geq 0$,
$\Phi_{L}(0) = 0$, $\Phi_{L}^{''} - L\Phi_{L}^{'} = 1 $, $\Phi_{L}^{'}(x) \ge$ 0 and $\Phi_{L}^{''}(x)  \ge 1$, if $x \ge$ 0. Since $V^{n,p} \cdot M$ and $N^{n,p}$ are square integrable martingales, their expectations are constant and equal to zero. Thanks to It\^o's formula applied to $\Phi_{L}(U^{n, p})$, we get
 \begin{tabbing}
 $  \mathbb{E}{\Phi_{L}(U_{0}^{n, p})}  -  \mathbb{E}{\Phi_{L}(U_{T}^{n, p}}) =$ \= $ \mathbb{E}\disp{\int_{0}^{T}\big(\Phi_{L}^{'}(U_{s}^{n, p})(g^{n}(s,U_{s}^{n},V_{s}^{n})- g^{p}(s,U_{s}^{p},V_{s}^{p}))\big)dC_{s}}$ \\
$ $\\
\> $   - \mathbb{E}\disp{\int_{0}^{T}\frac{\Phi_{L}^{''}}{2}(U_{s}^{n, p})|m_{s}(V_{s}^{n, p})|^{2}dC_{s} - \mathbb{E}\int_{0}^{T}\frac{\Phi_{L}^{''}}{2}(U_{s}^{n, p})d\lan{N}^{n, p}\ran_{s}} .  $\\
\end{tabbing}
 Then, since both $g^{n}$ and $g^{p}$ satisfy ($H_{1}^{'}$) with the same parameters, 
\begin{tabbing} 
\= $ \disp{ |g^{n}(s,U_{s}^{n},V_{s}^{n})}$ \= $ \disp{ - g^{p}(s,U_{s}^{p},V_{s}^{p})|  } $\\      
\> $ $ \\
 \> $  \disp{  \leq 2\bar{\alpha}_{s}+ \frac{\gamma}{2}|m_{s}(V_{s}^{n})|^{2} + \frac{\gamma}{2}|m_{s}(V_{s}^{p})|^{2} } $\\ 
\> $ $\\
 \> $ \disp{\leq 2\bar{\alpha}_{s} }$ \= $ \disp{ + \frac{3\gamma}{2}\big(|m_{s}(V_{s}^{n, p} )|^{2}  + |m_{s}(V_{s}^{p} - \tilde{V}_{s})|^{2} + |m_{s}\tilde{V}_{s}|^{2} \big) + \gamma \big(|m_{s}(V_{s}^{p} - \tilde{V}_{s})|^{2} + |m_{s}\tilde{V}_{s}|^{2}\big) }$\\
\> $ $\\
\> $ \disp{ \leq 2\bar{\alpha}_{s} + \frac{3\gamma}{2}\big(|m_{s}(V_{s}^{n, p} )|^{2}\big)  + \frac{5\gamma}{2}\big(|m_{s}(V_{s}^{p} - \tilde{V}_{s})|^{2} + |m_{s}\tilde{V}_{s}|^{2} \big)}. $\\
\end{tabbing}
The two last inequalities result from the convexity of: $z \to |z|^{2}$. Using these estimates and transferring $$\mathbb{E}\big(\disp{\int_{0}^{T}\frac{\Phi_{L}^{''}}{2}(U_{s}^{n, p})|m_{s}(V_{s}^{n, p})|^{2}dC_{s}}\big) \; \; \textrm{and} \; \; \mathbb{E}\big(\disp{\int_{0}^{T}\Phi_{L}^{'}(U_{s}^{n, p})\frac{3\gamma}{2}|m_{s}(V_{s}^{n, p})|^{2}dC_{s}}\big), $$ in the left-hand side of It\^o's formula applied to $\Phi_{L}(U^{n, p})$, we obtain 
\begin{tabbing}
\=   $  \disp{ \mathbb{E}{\Phi_{L}(U_{0}^{n, p})}} \;  +  \disp{\frac{1}{2}\mathbb{E}\big(|N_{T}^{n, p}|^{2}\big) +  \mathbb{E}\int_{0}^{T}( (\frac{\Phi_{L}^{''}}{2} - \frac{3\gamma}{2} \Phi_{L}^{'})(U_{s}^{n, p})|m_{s}(V_{s}^{n, p})|^{2}dC_{s})} $\\
$ $ \\
\> $ \quad \disp{  \leq  \; \mathbb{E}{\Phi_{L}(\tilde{B}^{n} - \tilde{B}^{p})} + \mathbb{E}\int_{0}^{T}\Phi_{L}^{'}(U_{s}^{n, p})\left(2\bar{\alpha}_{s} + \frac{5\gamma}{2}(|m_{s}(V_{s}^{p} - \tilde{V}_{s})|^{2} + |m_{s}\tilde{V}_{s}|^{2})\right) dC_{s}.} \;  (**)$    \\  
 \end{tabbing}
Setting: $L := 8\gamma $, and using the definition (\ref{eq: PhiL}), we claim \\
\begin{equation}\label{eq: conditionL}
\Phi_{L}^{''} - 8\gamma \Phi_{L}^{'} = 1,
\end{equation}
which entails the 
positiveness of the last term of the left-hand side. Then, thanks to
the weak convergence of ($V^{n}$) to $\tilde{V}$ (and of ($N^{n}$) to $\tilde{N}$) and the convexity of $z \to |z|^{2}$, we have\\
\begin{eqnarray}\label{eq:  convfaible}\displaystyle{\liminf_{n \to
      \infty}\mathbb{E}\int_{0}^{T}( (\frac{\Phi_{L}^{''}}{2}-  \frac{3\gamma}{2} 
    \Phi_{L}^{'})(U_{s}^{n, p})|m_{s}(V_{s}^{n, p})|^{2})  dC_{s} }  \geq   \\
\nonumber  \mathbb{E}\int_{0}^{T}( (\frac{\Phi_{L}^{''}}{2}  -   \frac{3\gamma}{2} 
  \Phi_{L}^{'}) ( \tilde{U}_{s} - U_{s}^{p})( |m_{s}(\tilde{V}_{s} - V_{s}^{p})|^{2})dC_{s}).  \\ 
\nonumber \end{eqnarray}
 Similarly, we get
\begin{eqnarray}\label{eq:  convfaible2}  \displaystyle{\liminf_{n \to
      \infty}\mathbb{E}\big(|N_{T}^{n, p}|^{2}\big)  }  \geq   \mathbb{E}\big(|\tilde{N}_{T} - N_{T}^{p}|^{2}\big)  .  \\ 
\nonumber \end{eqnarray}
Using the almost sure convergence of the increasing sequence ($U^{n}$) to $\tilde{U}$, the dominated convergence theorem yields
\begin{eqnarray}\label{eq: PALimit}
\nonumber \disp{\Phi_{L}^{'}(U_{s}^{n, p})\big(\frac{5\gamma}{2}(|m_{s}(\tilde{V}_{s} - V_{s}^{p})|^{2} + |m_{s}\tilde{V}_{s}|^{2}) + 2\bar{\alpha}_{s}\big)}\quad  \quad \quad \quad \quad \quad \quad  &  \\
\quad \quad \quad \disp{\leq \; \Phi_{L}^{'}(\tilde{U}_{s}-  U_{s}^{p})\big(\frac{5\gamma}{2}(|m_{s}(\tilde{V}_{s} - V_{s}^{p})|^{2} + |m_{s}\tilde{V}_{s}|^{2}) + 2\bar{\alpha}_{s}\big),}&\\ 
\nonumber  \end{eqnarray}  which holds uniformly in $n$. Besides, the process in the right-hand side of (\ref{eq: PALimit}) is integrable w.r.t. $dC$,
as a product of a bounded process and a sum of integrable processes.
 Then, we use both (\ref{eq: convfaible}) and (\ref{eq: convfaible2}) to give a lower bound of the left-hand side of inequality (**). For the right-hand side of (**), we rely on (\ref{eq: PALimit}) and on the almost sure and increasing convergence of $(\tilde{B}^{n})$ to $\tilde{B}$ to get
\begin{tabbing}
  \= $ \disp{ \mathbb{E}\Phi_{L}(\tilde{U}_{0} - U_{0}^{p})}  + \disp{ \frac{1}{2}\mathbb{E}
 \big(|\tilde{N}_{T} - N_{T}^{p}|^{2}\big)} $\\
\> $  \; \disp{ + \; \mathbb{E}\int_{0}^{T}( (\frac{\Phi_{L}^{''}}{2} - \frac{3\gamma}{2}\Phi_{L}^{'})(\tilde{U}_{s}- U_{s}^{p})|m_{s}(\tilde{V}_{s} - V_{s}^{p})|^{2}dC_{s})}  $ \\
\> $\quad  \disp{ \leq \; \mathbb{E}\left( \Phi_{L}(\tilde{B} - \tilde{B}^{p}) +  \int_{0}^{T} \Phi_{L}^{'}( \tilde{U}_{s}-
   U_{s}^{p} )(\frac{5\gamma}{2}|m_{s}(\tilde{V}_{s} - V_{s}^{p})|^{2} + 2\bar{\alpha_{s}} + \frac{5\gamma}{2}|m_{s}\tilde{V}_{s}|^{2}) dC_{s}\right)}. $\\
\end{tabbing}
 Transferring now $\mathbb{E}\big(\int_{0}^{T} \Phi_{L}^{'}( \tilde{U}_{s}-
   U_{s}^{p} )(\frac{5\gamma}{2}|m_{s}(\tilde{V}_{s} - V_{s}^{p})|^{2}) dC_{s}\big)$ in the left-hand side of this inequality
and using properties of $\Phi_{L}$ and, in particular, (\ref{eq: conditionL}), we obtain
\[ \begin{array}{l}
\disp{ \mathbb{E}\Phi_{L}(\tilde{U}_{0} - U_{0}^{p}) + \frac{1}{2}\mathbb{E}\left(\int_{0}^{T}|m_{s}(\tilde{V}_{s} -
   V_{s}^{p})|^{2}dC_{s} + |\tilde{N}_{T} - N_{T}^{p}|^{2}\right) }\\

\quad \; \leq  \disp{\mathbb{E}\left( \Phi_{L}(\tilde{B} - \tilde{B}^{p}) +  \int_{0}^{T} \Phi_{L}^{'}( \tilde{U}_{s}-
   U_{s}^{p} )(2\bar{\alpha_{s}} + \frac{5\gamma}{2}|m_{s}\tilde{V}_{s}|^{2}) dC_{s}\right)}.\\
 \end{array} \]
  Thanks to the convergence of $(\tilde{U}_{s} - U^{p}_{s})$ to $0$ (holding true $\mb{P}$-a.s. and for all $s$) and since $|m \tilde{V}|^{2}$ and $\bar{\alpha}$ are in $L^{1}(dC \otimes d\mb{P})$, the dominated convergence theorem  entails the convergence of the right-hand side to $0$. Taking the limit sup over $p$ in the left-hand side, it yields
$$ \disp{\lim \sup_{p \to \infty} \mb{E}\left(\frac{1}{2}\mathbb{E}\left(\int_{0}^{T}|m_{s}(\tilde{V}_{s} -
   V_{s}^{p})|^{2}dC_{s} + |\tilde{N}_{T} - N_{T}^{p}|^{2}\right) \right)} \le 0,$$
which ends the proof.\\
\section{Applications to finance }
In this section, we study the problem (\ref{eq: utilmaxim}) stated in the introduction for three types of utility functions.
\subsection{The case of the exponential utility }
 \begin{theorem}\label{dynammethod} 
$\bullet \;$ For any fixed $t$, the value function $x \to V_{t}^{B}(x)$ can be expressed in term of the unique solution ($Y,Z,L$) of BSDE of type $(\textrm{Eq1})$ given by ($F^{\alpha}, \beta, B$) 
\begin{equation}\label{eq: fonctionval} V_{t}^{B}(x)= U_{\alpha}(x - Y_{t}).\end{equation}
 $ \beta:= \alpha$ corresponds to the risk-aversion parameter, $B$ is the contingent claim and $F^{\alpha}$ is the generator, whose expression is 
\begin{eqnarray*}
 F^{\alpha}(s, z) = \displaystyle{  \mathop{\inf_{\nu \in \mathcal{C}}\big(\frac{\alpha}{2}|m_{s}(\nu - (z +  \frac{\lambda_{s}}{\alpha}))|^{2}\big) }}
- (m_{s}z)^{'}(m_{s}\lambda_{s})  - \frac{1}{2\alpha}|m_{s}\lambda_{s}|^{2}.  \\
\nonumber  \end{eqnarray*}
 $\bullet \;$ There exists an optimal strategy $\nu^{*}:= (\nu_{s}^{*})_{s \in [t, \; T]}$ such that: $\nu^{*} \in \mc{A}_{t}$, and satisfying, $\mathbb{P}\textrm{-a.s.} $ and for all $s$, 
\begin{equation}\label{eq: optimstrat} \nu_{s}^{*} \in  \textrm{arg}\displaystyle{\min_{ \nu \in \mathcal{C}} |m_{s}(\nu - (Z_{s} + \frac{\lambda_{s}}{\alpha}))|^{2}}.
 \end{equation}
 $\bullet \;$ Extending the definition of $V_{t}^{B}(x)$ to an arbitrary stopping time $\tau$, we set
$$ V_{\tau}^{B}(x) : = \textrm{ess} \disp{\sup_{\nu} \mb{E}^{\mc{F}_{\tau}}\left(U_{\alpha}(x +\int_{\tau}^{T}\disp{\sum_{i} \nu_{u}^{i}\frac{dS_{u}^{i}}{S_{u}^{i}}} - B) \right)},$$
where, in this expression, trading strategies $\nu$ are defined on $[\tau, T]$. Then, for any $\tau$,
$$ V_{\tau}^{B}(x) = U_{\alpha}(x -Y_{\tau}) = R_{\tau}^{\nu^{*}},$$
 and we recover the formulation of the dynamic programming principle for $V^{B}$
\begin{equation}\label{eq: Pdynam} \forall \; \tau,\; \sigma,  \tau \leq \sigma, \;\mc{F}\textrm{-stopping time}, \quad V_{\tau}^{B}(x) = \mb{E}^{\mc{F}_{\tau}}(V_{\sigma}^{B}(X_{\sigma}^{\nu^{*}, \tau, x})).
\end{equation}
\end{theorem}
\textbf{Remark:}
 To give sense to the expression $V_{\sigma}^{B}(X_{\sigma}^{\nu^{*}, \tau, x}) $, we refer to the footnote given at the bottom of page 5: in fact, $X_{\sigma}^{\nu^{*}, \tau, x}:= x + \disp{\int_{\tau}^{\sigma}\nu_{u}^{*}\frac{dS_{u}}{S_{u}}}$, is an attainable wealth at time $\sigma$, when starting from $x$ at time $\tau$. \\
\textbf{Proof:}
To prove (\ref{eq: fonctionval}), we first rely on the results obtained in Section 3 to claim the existence of a unique solution ($Y, Z, U$) of the BSDE($F^{\alpha}, \alpha, B$). Then, we give the expression of $R^{\nu} := U_{\alpha}(X^{\nu} - Y)$ obtained in the last paragraph of Section 2.1
 $$\forall \; s \in [t, T], \quad R_{s}^{\nu}  = R_{t}^{\nu}\tilde{M}_{t, s}^{\nu}\exp(A_{s}^{\nu} - A_{t}^{\nu}),$$
with: $ \tilde{M}_{t, s}^{\nu} := \mc{E}_{t, s}(-\alpha(\nu - Z) \cdot M + \alpha L).$
Since the continuous stochastic exponential is a positive local martingale and since $A^{\nu} \geq 0$, there exists a sequence of stopping time 
($\tau_{n}$) such that ($R_{\cdot \wedge \tau_{n}}^{\nu}$) is a supermartingale (for each $\nu$), which entails \\ \[\forall \; s, \;t \leq s \leq T,\; \forall \;A \in \mathcal{F}_{t}, \quad \mathbb{E}\big(R^{\nu}_{s \wedge \tau_{n}}\mathbf{1}_{A}\big) \leq \mathbb{E}\big(R^{\nu}_{t \wedge \tau_{n}}\mathbf{1}_{A}\big) .\] 
 Using the definition of admissibility and the boundedness of $Y$, we obtain the uniform integrability of ($R^{\nu}_{t \wedge \tau_{n}}$) and ($R^{\nu}_{s \wedge \tau_{n}}$). Passing to the limit, we get: $\mathbb{E}\big(R^{\nu}_{s}\mathbf{1}_{A}\big) \leq \mathbb{E}\big(R^{\nu}_{t}\mathbf{1}_{A}\big)$, which entails the supermartingale property of $R^{\nu}$, as soon as: $\nu \in \mc{A}_{t}$. Both this supermartingale property and the relation: $R_{t}^{\nu} = U_{\alpha}(x- Y_{t})$, entail 
$$ V_{t}^{B}(x) = \textrm{ess} \disp{\sup_{\nu \in \mc{A}_{t}} \mb{E}^{\mc{F}_{t}}(U_{\alpha}(X_{T}^{\nu, x,t} -B)}
\le U_{\alpha}(x- Y_{t}).
 $$
\indent Now, to obtain the equality
 (\ref{eq: fonctionval}), we focus on
 the second point of Theorem \ref{dynammethod}. Firstly, the infimum in the expression of $F^{\alpha}$ exists, since: $z \to F^{\alpha}(s, z)$ is a continuous functional of $z$, which tends to $+\infty$, as $|z|$ goes to $\infty$. Furthermore, relying on the same selection argument as in lemma 11 in \cite{ImkelleretHu} and thanks to the continuity of the functional and the predictability of the processes $\lambda$ and $Z$, there exists a measurable choice of $\nu_{s}^{*}$ satisfying (\ref{eq: optimstrat}), i.e. $A^{\nu^{*}} \equiv 0$ . Finally, to check that: $\nu^{*} \in \mc{A}_{t}$, we first argue that, from the choice of $\nu^{*}$ given in Theorem \ref{dynammethod} and since 0 is in $\mathcal{C}$, 
$$\forall \; s \in [0, \; T], \quad |m_{s}(\nu_{s}^{*} - (Z_{s} + \frac{\lambda_{s}}{\alpha}))| \leq |m_{s}(Z_{s} + \frac{\lambda_{s}}{\alpha})|.$$
Noting that: $ |m_{s}(\nu_{s}^{*} - Z_{s}))| \leq |m_{s}(\nu_{s}^{*} - (Z_{s} + \frac{\lambda_{s}}{\alpha}))| + |m_{s}\frac{\lambda_{s}}{\alpha}|$,
we obtain a control of $|m(\nu^{*} - Z))|$ depending only on the processes $Z$ and $\lambda$. Hence and thanks to Kazamaki's criterion (see \cite{Kazamaki}), $\mathcal{E}(-\alpha(\nu^{*} - Z) \cdot M)$ is a true martingale.
The process $R^{\nu^{*}}$, such that, for all $s$,
 $$R_{s}^{\nu^{*}} = -e^{-\alpha(x - Y_{t})}\mathcal{E}_{t, s}(-\alpha(\nu^{*} - Z) \cdot M + \alpha L),$$ is a true martingale, which implies that: $\nu^{*} \in \mathcal{A}_{t}$ and the equality (\ref{eq: fonctionval})
.\\
\indent To recover the dynamic principle, we define the $\mc{F}_{\tau}$-measurable random variable $V_{\tau}^{B}(x)$  the same way as $V_{t}^{B}(x)$ and for any $\mc{F}$-stopping time $\tau$. The same procedure as the one used to prove (\ref{eq: fonctionval}) entails
 $$V_{\tau}^{B}(x) = U_{\alpha}(x - Y_{\tau}) = U_{\alpha}(X_{\tau}^{\nu^{*}, \tau, x} - Y_{\tau}) = R_{\tau}^{\nu^{*}}.$$  Applying the optional sampling theorem between $\tau$ and $\sigma $ to the martingale $R^{\nu^{*}} :=U_{\alpha}(X^{\nu^{*}, \tau, x} - Y)$, we get (\ref{eq: Pdynam}).
\begin{flushright}
 $\square$
\end{flushright}

\subsection{Power and logarithmic utilities}
\hspace{0.5cm} As in \cite{ImkelleretHu}, we introduce two other types of utility functions:\\
\begin{itemize}
\item The first one is the power utility, defined for all real $\gamma$, $\gamma \in ]0, 1[$, by:
$U_{\gamma}(x) = \frac{1}{\gamma} x^{\gamma}$ ($\gamma$ being fixed, we will write $U^{1}$ instead of $U_{\gamma}$).\\
\item The second one is the logarithmic utility, given by: 
 $U^{2}(x) = \ln(x)$.\\
\end{itemize}
  Contrary to the exponential case, we have to impose that the wealth process is positive. Besides, in these two cases, there is no liability any more (i.e. $B \equiv 0$, in the problem (\ref{eq: utilmaxim})). We provide here another notion of strategy: a constrained trading strategy is a
$d$-dimensional process $\rho$ taking its values in the constraint set $\mc{C}$ and such that each component $\rho_{i}$ stands for the part of the wealth
invested in stock $i$. The discounted price process $S$ is again assumed to satisfy (\ref{eq: SDE1}) and we denote by: $X^{\rho}:= X^{\rho,t,x}$, the wealth process associated with strategy $\rho$ and such that: $X_{t}^{\rho} = x$. Its expression for any $s$, $s\in [t, T]$, is
$$ X_{s}^{\rho} = x + \int_{t}^{s}X_{u}^{\rho}\rho_{u}\frac{dS_{u}}{S_{u}} = x +
\int_{t}^{s}X_{u}^{\rho}\rho_{u}dM_{u} + \int_{t}^{s} X_{u}^{\rho}\rho_{u}^{'}d\lan{M} \ran_{u}\lambda_{u}.
$$
 For each case, we will give a definition of the admissibility set for trading strategies (this set is always denoted by $\mc{A}_{t}$).
Denoting by $U$ the utility function, we are going to characterize the value function $x \to V_{t}(x)$ at time $t$, which is defined by
\begin{equation}\label{eq: pboptim} V_{t}(x) = \textrm{ess} \displaystyle{\sup_{\rho, \;\rho \in \mc{A}_{t} }\mathbb{E}^{\mathcal{F}_{t}}(U(x + \int_{t}^{T}X_{u}^{\rho}\rho_{u}\frac{dS_{u}}{S_{u}}))}. \end{equation}

\subsubsection{The power utility case} 
\begin{definition}
The set of admissible strategies $\mc{A}_{t}$ consists of all $d$-dimensional predictable processes $\rho:=(\rho_{s})_{s \in [t, T]}$ such that $\rho_{s} \in \mc{C}$ ($\mb{P}$-a.s. and for all $s$) as well as
$$\int_{t}^{T}\rho_{s}^{'}d\lan{M} \ran_{s}\rho_{s} = \int_{t}^{T}|m_{s}\rho_{s}|^{2}dC_{s} < \infty, \; \mathbb{P} \textrm{-}\;\textrm{a.s}. $$
\end{definition}
\indent This condition entails that the stochastic exponential $\mathcal{E}(\rho \cdot M)$ is a continuous local martingale. We can now solve the problem (\ref{eq: pboptim}) for the power utility function $U^{1}$.
\begin{theorem}\label{caspuissance}
Let $x \to V_{t}^{1}(x)$ be the value function associated with the problem (\ref{eq: pboptim}) and with $U = U^{1}$.\\
$ \bullet \;$ Its expression is 
\[ V^{1}_{t}(x) = \frac{x^{\gamma}}{\gamma} \exp(Y_{t}),  \]
where ($Y, Z, L$) stands for the unique solution of the BSDE($f^{1}, \frac{1}{2}, 0$) of type ($\textrm{Eq1}$)  
\begin{eqnarray*}
 Y_{t} = 0 -\int_{t}^{T}f^{1}(s, Z_{s})dC_{s} + \int_{t}^{T}\frac{1}{2}d\lan{L} \ran_{s} 
- \int_{t}^{T}Z_{s}dM_{s} - (L_{T} - L_{t}), \\ \end{eqnarray*}
and where $L$ is a real martingale strongly orthogonal to $M$. The expression of
$f^{1}$ is
  \begin{eqnarray}\label{eq: generator1}
\nonumber f^{1}(s, z) =    \displaystyle{\inf_{\rho,\; \rho \in \mc{C}}  \frac{\gamma(1 - \gamma )}{2}\left( |m_{s}(\rho - (\frac{z + \lambda_{s}}{1 - \gamma }))|^{2}  \right) }  \\
   -\frac{\gamma(1 - \gamma )}{2}|m_{s}(\frac{z + \lambda_{s}}{1 - \gamma })|^{2} 
  - \frac{1}{2}|m_{s}z|^{2}.  \\
\nonumber \end{eqnarray} 
$\bullet \;$ There exists an optimal strategy $\rho_{1}^{*}$ satisfying, $\mathbb{P}$-a.s. and for all $s$, \\ \begin{equation}\label{eq: optimstratpuiss} \disp{(\rho_{1}^{*})(s) \in \textrm{arg}\displaystyle{\min_{ \rho,\; \rho \in \mc{C}}|m_{s}(\rho - (\frac{Z_{s} + \lambda_{s}}{1 - \gamma }))|^{2} }}.
  \end{equation}
\end{theorem}
\textbf{Remark:}
 The expression of the optimal strategy $\rho^{*}$ is already known in the brownian setting and when there is no trading constraints: for instance, we refer the reader to the expression (3.19) given in \cite{Zhariphopoulou}.
In this paper, the wealth process $X^{\pi}$ satisfies
\begin{equation}\label{eq: wealthp}
 dX_{s}^{\pi}:= r X_{s}^{\pi}ds + X_{s}^{\pi}\big(\sigma_{s}\pi_{s} dW_{s} + (\mu - r)\pi_{s}ds\big).
\end{equation}
To obtain (3.19), we just take the correlation factor (denoted by $\rho$ in \cite{Zhariphopoulou}) equal to $0$ and we replace $\rho^{*}$ by $\frac{\pi^{*}\sigma}{x}$ (here, $\rho^{*}$ stands for the proportion invested in the risky asset, whereas, in \cite{Zhariphopoulou}, $\pi^{*}$ stands for the amount of wealth), and $\lambda$ by $\sigma^{-1} (\mu - r)$. 
In the case of constant coefficients in (\ref{eq: wealthp}), we recover that the optimal proportion is equal to $\frac{\mu - r}{(1 - \gamma)\sigma^{2}}$.\\
\textbf{Proof:}
We just give here the sketch of the proof, which is similar to the one given in the exponential case and relies on the same dynamic method as in \cite{ImkelleretHu}. To this end, we define the process $R^{\rho}$ for all $s$, $s \in [t, T]$, by: $R_{s}^{\rho} = X_{s}^{\rho}\exp(Y_{s})$.
We first write
\[ X_{s}^{\rho} = x + \int_{t}^{s}X_{u}^{\rho}\rho_{u}dM_{u} + \int_{t}^{s}X_{u}^{\rho}(m_{u}\rho_{u})^{'}(m_{u}\lambda_{u})dC_{u}, \]
and since $Y$ is solution of the BSDE($f^{1}, \frac{1}{2}, 0$), it results from simple computations
$$ R_{s}^{\rho} = R_{t}^{\rho}\frac{1}{\gamma}\mathcal{E}_{t, s}((\gamma\rho + Z) \cdot M + L)\exp(\tilde{A}_{s}^{\rho} - \tilde{A}_{t}^{\rho}),  $$
where the process $\tilde{A}^{\rho}$ is such that  \\
 $$\tilde{A}_{s}^{\rho} = \int_{0}^{s} \big( f^{1}(u, Z_{u}) + \frac{1}{2}|m_{u}Z_{u}|^{2}  + \frac{\gamma(\gamma -1)}{2}|m_{u}\rho_{u}|^{2} + \gamma(m_{u}\rho_{u})^{'}(m_{u}(Z_{u} + \lambda_{u})) \big)dC_{u}. $$ 
By the definition of $f^{1}$, we check:\\
\textbullet$\;$ $R^{\rho}$ is a supermartingale for any $\rho$, $\rho \in \mc{A}_{t}$,\\
\textbullet$\;$ $R^{\rho^{*}}$ is a martingale for any strategy $\rho_{1}^{*}$ satisfying (\ref{eq: optimstratpuiss}), taking into consideration that, for such a strategy, we have: $|m \rho_{1}^{*}| \leq |m \frac{(Z + \lambda)}{(1 - \gamma)}|$. \\
Besides, we obtain 
$$V^{1}_{t}(x) =  \mb{E}^{\mc{F}_{t}}(R_{T}^{\rho_{1}^{*}}) = R_{t}^{\rho_{1}^{*}} = \frac{x^{\gamma}}{\gamma} \exp(Y_{t}). $$

\subsubsection{The logarithmic utility case} 
 Once again, we introduce the notion of admissible strategy adapted to our problem.
\begin{definition}
The set of admissible strategies $\mc{A}_{t}$ consists of all $d$-dimensional predictable processes $\rho$ such that, $\rho_{s} \in \mc{C}$, $\mb{P}$-a.s. and for all $s$, and such that $$\disp{ \mathbb{E}\big(\int_{t}^{T}\rho_{s}^{'}d\lan{M} \ran_{s}\rho_{s}\big) = \mathbb{E}\big(\int_{t}^{T}|m_{s}\rho_{s}|^{2}dC_{s}\big) < \infty} .$$
\end{definition}
\begin{theorem}\label{caslog}
 Let $x \to V_{t}^{2}(x)$ be the value function associated with the problem (\ref{eq: pboptim}) and with $U = U^{2}$ for utility function.\\
$ \bullet \; $ Its expression is  $V^{2}_{t}(x) := \ln(x) + Y_{t}, $
where $Y$ stands for the unique solution of the BSDE($f^{2}, 0$) of type ($\textrm{Eq2}$) 
\[  Y_{t} = 0 - \int_{t}^{T}f^{2}(s)dC_{s} 
- \int_{t}^{T}Z_{s}dM_{s} - \int_{t}^{T}dL_{s}, \]
and where the expression of $f^{2}$ is
\begin{eqnarray}\label{eq: generator2}
f^{2}(s) =  \displaystyle{\inf_{\rho,\; \rho \in \mc{C}} \frac{1}{2}|m_{s}(\rho - \lambda_{s})|^{2}} - \frac{1}{2}|m_{s}\lambda_{s}|^{2}. \\
\nonumber \end{eqnarray}
$ \bullet \; $ There exists an optimal strategy $\rho_{2}^{*}$ satisfying ($\mathbb{P}$-a.s. and for all $s$) \begin{equation}\label{eq: optimstrat2} (\rho_{2}^{*})(s) \in \textrm{arg}\displaystyle{\min_{ \rho, \; \rho \in \mathcal{C} }|m_{s}(\rho -\lambda_{s})|^{2} }.
\end{equation}
\textbf{Remark:}
 As in the power utility case, we recover the expression of the optimal proportion in the brownian setting. Assuming that the coefficients $\mu$, $\sigma $ and $r$ are constants, this proportion is equal to: $\rho^{*} \equiv \frac{(\mu -r)}{\sigma^{2}}$.

\end{theorem}
\textbf{Proof:} $\hspace{0.5cm}$
 The wealth process $X^{\rho}$ satisfies again
\[ X_{s}^{\rho} = x + \int_{t}^{s}X_{u}^{\rho}\rho_{u}dM_{u} + \int_{t}^{s}X_{u}^{\rho}(m_{u}\rho_{u})^{'}(m_{u}\lambda_{u})dC_{u}. \]
Now, It\^o's formula and the assumption that $Y$ solves a BSDE of type ($\textrm{Eq2}$) yield  
\[ R_{s}^{\rho} = \ln(X_{s}^{\rho}) + Y_{s} = \ln(x) + Y_{t} + \int_{t}^{s}\big( (\rho_{u} + Z_{u})dM_{u} + dL_{u}\big) + A_{2}^{\rho}(s) - A_{2}^{\rho}(t), \]
where the process $A_{2}^{\rho}$ is given by $$ A_{2}^{\rho}(s) = \int_{0}^{s}(f^{2}(u) - \frac{1}{2}|m_{u}\rho_{u}|^{2}+ (m_{u}\rho_{u})^{'}(m_{u}\lambda_{u}))dC_{u}.$$
From the definition of $f^{2}$, we obtain: 
$A_{2}^{\rho} \leq 0$, and we deduce:\\
$\bullet \;$ $\ln(X^{\rho}) + Y$ is a supermartingale, for any $\rho$ such that $\rho \in \mc{A}_{t} $.\\
 If, besides, $\rho_{2}^{*}$ satisfies (\ref{eq: optimstrat2}) then: $A_{2}^{\rho_{2}^{*}} = 0 $ and, hence: $|m (\rho_{2}^{*}- \lambda)| \leq |m \lambda|$. The assumption ($H_{\lambda}$) on $\lambda$ implies the uniform integrability of $R^{\rho_{2}^{*}}$, and it follows that  \\
$\bullet \; $  $\ln(X^{\rho_{2}^{*}}) + Y$ is a martingale.\\
Hence, such a strategy $\rho_{2}^{*} $ is optimal and
applying the optional sampling theorem to $R^{\rho_{2}^{*}}$, it implies
$$ V^{2}_{t}(x) = \mb{E}^{\mc{F}_{t}}(R_{T}^{\rho_{2}^{*}}) = R_{t}^{\rho_{2}^{*}} = \ln(x) + Y_{t}. $$
\begin{flushright}
 $\square$
\end{flushright} 
 
\section{Conclusion}
\hspace{0.5cm} In this paper, we have solved the utility maximization problem by computing the value function and characterizing all optimal strategies: the novelty of our study is that we have used a dynamic method in the context of a general (and non necessarily Brownian) filtration and in presence of portfolio constraints. This last assumption entails that the introduced BSDEs have quadratic growth. \\
\indent Since we are not in the Brownian setting, the first part of our work consists in justifying new existence and uniqueness results for solutions of a type of quadratic BSDEs. This study leads to an expression of the value function in terms of a solution of a BSDE of the previous type. Relying on the dynamic principle, we are able to characterize the value function for three cases of utility functions. This type of BSDE has already been studied in a particular case in \cite{ManiaetSchw} in connection with the notion of the indifference utility price. However, one of the main difference in \cite{ManiaetSchw} is that no constraints are imposed on the portfolio. Furthermore and contrary to our setting, they refer to duality methods. Our study depends heavily on the assumption that the filtration is continuous and we hope to study the case when jumps are allowed. Another perspective is to study the connection with the problem of utility indifference pricing.

\section{Appendix: proof of proposition \ref{estimclass}}
Contrary to lemma \ref{estimapriori}, where the process $Y$ is supposed to be in $S^{\infty}$, in this proposition, the process $U^{n}$ is only assumed to be in $S^{2}$. We first apply It\^o's formula for $(e^{\Gamma C_{t}} |U_{t}^{n}|^{2})$, 
 $\Gamma$ being a non negative constant
  \begin{equation}\label{eq: formuleIto1}
d(e^{\Gamma C_{t}} |U_{t}^{n}|^{2}) = \Gamma e^{\Gamma C_{t}} |U_{t}^{n}|^{2}dC_{t} + e^{\Gamma C_{t}}\big(2U_{t}^{n} dU_{t}^{n} + d\lan{U^{n}}\ran_{t} \big), 
  \end{equation}
with 
 \begin{tabbing}
 $ 2U_{t}^{n} dU_{t}^{n} + d\lan{U^{n}}\ran_{t} $\=$:= -2U_{t}^{n}g^{n}(t, U_{t}^{n}, V_{t}^{n})dC_{t}  + |m_{t}V_{t}^{n}|^{2}dC_{t} + d\lan{N^{n}} \ran_{t}$\\
 \\
$ $ \> $\; \quad \;+  \;2U_{t}^{n} \big(V_{t}^{n}dM_{t} + dN_{t}^{n} \big). $ \\
 \end{tabbing}
Since ($U^{n}, V^{n}, N^{n}$) is in $S^{2} \times L^{2}(d \lan{M} \ran \times d\mathbb{P}) \times \mathcal{M}^{2}([0,
T])$, we can prove that the following process \\
\begin{equation}\label{eq: formeK}
\forall \; s \in  [0,T], \quad K_{s}:= \disp{\int_{0}^{s}2 e^{\Gamma C_{u}}U_{u}^{n} \big(V_{u}^{n}dM_{u} + dN_{u}^{n} \big)}, 
\end{equation}
 is a true martingale.
 We now fix $t$ ($t \in [0,  T]$) 
and we rewrite It\^o's formula (\ref{eq: formuleIto1}) in its integrated form between $s$ ($t \le s \le T$) and $T$ 
  \[ \begin{array}{ll}
 e^{\Gamma C_{s}} |U_{s}^{n}|^{2} -e^{\Gamma C_{T}} |U_{T}^{n}|^{2} &
 =\;\disp{\int_{s}^{T}e^{\Gamma C_{u}} U_{u}^{n}\big(-\Gamma U_{u}^{n} + \;2g^{n}(u, U_{u}^{n}, V_{u}^{n})\big)dC_{u} }\\
\\
 &- \disp{\int_{s}^{T}e^{\Gamma C_{u}}\big(|m_{u}V_{u}^{n}|^{2}dC_{u} +d\lan{N^{n}\ran}_{u} \big) -\big(K_{T} - K_{s} \big). } \\
  \end{array} \]
 We rely on the Lipschitz property of the generator $g^{n}$ to get 
 $$ 2|U_{u}^{n}||g^{n}(u, U_{u}^{n}, V_{u}^{n} )| \leq 2|U_{u}^{n}||g^{n}(u, 0, 0)| + 2L_{n}\big( |U_{u}^{n}|^{2} + |U_{u}^{n}||m_{u}V_{u}^{n} |\big),$$
 and using the inequality: $|2L_{n}ab| \le ( 2(L_{n})^{2}a^{2} + \frac{1}{2}b^{2})$, we obtain
 $$ 2L_{n}|U_{u}^{n}||m_{u}V_{u}^{n} | \le 2(L_{n})^{2}|U_{u}^{n}|^{2} + \frac{1}{2}|m_{u}V_{u}^{n}|^{2}. $$
Combining these two last inequalities, setting: $\Gamma = 2((L_{n})^{2} + L_{n})$, and taking the expectation w.r.t $\mc{F}_{t}$ in It\^o's formula applied to $ e^{\Gamma C_{s}} |U_{s}^{n}|^{2}$ between $t$ and $T$, we get
 \begin{tabbing}
$\disp{ e^{\Gamma C_{t}} |U_{t}^{n}|^{2}  \le} $  \= $ \disp{\mb{E}\left(e^{\Gamma C_{T}} |U_{T}^{n}|^{2} |\mc{F}_{t}\right) }$\\
\\
$ $\> $\quad  \;\; +\; \disp{\mb{E}\left(\int_{t}^{T}e^{\Gamma C_{u}}\big(2|U_{u}^{n}||g^{n}(u, 0, 0)|  + \frac{1}{2}(|m_{u}V_{u}^{n}|^{2}) \big) dC_{u}|\mc{F}_{t}\right)} $ \\
\\
$ $\> $\quad \quad \;  - \disp{\mb{E}\left(\int_{t}^{T}e^{\Gamma C_{u}}\big(|m_{u}V_{u}^{n}|^{2}dC_{u} +d\lan{N^{n}\ran}_{u} \big)|\mc{F}_{t}\right) }.$ \\
\end{tabbing}
 This leads to \\
\newline
$\disp{\mb{E}\big(\int_{t}^{T}e^{\Gamma C_{u}}\big(|m_{u}V_{u}^{n}|^{2} dC_{u} + \lan{N} \ran_{u} \big)| \mc{F}_{t}\big)}$\\
\begin{equation}\label{eq: Ito2} \quad \quad\leq \disp{2\left(\mb{E}\big(e^{\Gamma C_{T}}|U_{T}^{n}|^{2} + 2\int_{t}^{T} e^{\Gamma C_{u}}|U_{u}^{n}||g^{n}(u, 0, 0)| dC_{u}| \mc{F}_{t}\big)\right).}
\end{equation}
We come back to It\^o's formula (\ref{eq: formuleIto1}) for the process $ e^{\Gamma C_{\cdot}} |U_{\cdot}^{n}|^{2}$ between $s$ and $T$. Taking then the supremum over $s$ ($s \in [t, T]$), it follows 
\[ \begin{array}{l}
 \disp{\sup_{ t\le s \le T}e^{\Gamma C_{s}} |U_{s}^{n}|^{2}} \le    \disp{e^{\Gamma C_{T}} |U_{T}^{n}|^{2}}
\\ \quad  \;+ \;\disp{2\int_{t}^{T}e^{\Gamma C_{u}}|U_{u}^{n}||g^{n}(u, 0, 0)| dC_{u}} 
\; + \disp{\sup_{ t\le s\le T}|K_{T} - K_{s}|} .\\
\end{array} \]
Applying the Burkholder-Davis-Gundy inequality to the supremum of the square integrable martingale $K$ and the relation: $ Cab \le \frac{C^{2}}{2}a^{2} +\frac{1}{2}b^{2}$, we deduce
the existence of a constant $C$ such that
\begin{tabbing}
$ $ \= $\mb{E}\left(\disp{\sup_{ t\le s \le T}e^{\Gamma C_{s}} |U_{s}^{n}|^{2}} |\mc{F}_{t} \right) \le \mb{E}\left( e^{\Gamma C_{T}} |U_{T}^{n}|^{2} + 2 \disp{\int_{t}^{T}e^{\Gamma C_{u}}|U_{u}^{n}||g^{n}(u, 0, 0)| dC_{u}|\mc{F}_{t} }\right)$\\
\\
$ $\>  $ \quad \quad + \frac{C^{2}}{2}\disp{\mb{E}\big(\int_{t}^{T}e^{\Gamma C_{u}}\big( |m_{u}V_{u}^{n}|^{2}dC_{u} + d\lan{N} \ran_{u} \big) |\mc{F}_{t}\big)  } + \frac{1}{2}\disp{\mb{E}\big(\disp{\sup_{ t\le s \le T}e^{\Gamma C_{s}} |U_{s}^{n}|^{2}} |\mc{F}_{t}\big)  }.$\\
\end{tabbing}
This constant $C$ is generic and may vary from line to line. 
Combining this last inequality with (\ref{eq: Ito2}), we deduce 
\begin{eqnarray*}
 \mb{E}\left(\disp{\sup_{ t\le s\le T}e^{\Gamma C_{s}} |U_{s}^{n}|^{2}} +\disp{ \int_{t}^{T}e^{\Gamma C_{u}}\big( |m_{u}V_{u}^{n}|^{2}dC_{u} + d\lan{N} \ran_{u} \big)}|\mc{F}_{t} \right) \\
 \quad \quad \quad \le C \mb{E}\left(e^{\Gamma C_{T}} |U_{T}^{n}|^{2} +  \disp{\int_{t}^{T}e^{\Gamma C_{u}}|U_{u}^{n}||g^{n}(u, 0, 0)| dC_{u}|\mc{F}_{t} }\right).\\
\end{eqnarray*}
To obtain the desired relation, we use a last estimate of the last term in the right-hand side of the previous inequality
\begin{tabbing}
$ $\=  $C\mb{E}\left(\disp{\int_{t}^{T}e^{\Gamma C_{u}}|U_{u}^{n}||g^{n}(u, 0, 0)|dC_{u}|\mc{F}_{t}} \right)$ \\   \\
\> $ \quad \;\le \; \frac{1}{2}\mb{E}\left(\disp{\sup_{ t\le u\le T}e^{\Gamma C_{u}} |U_{u}^{n}|^{2}|\mc{F}_{t}} \right) + \frac{C^{2}}{2}\mb{E}\left(\disp{\big(\int_{t}^{T}e^{\frac{\Gamma}{2} C_{u}}|g^{n}(u, 0, 0)|dC_{u}\big)^{2}|\mc{F}_{t}} \right).$\\
\end{tabbing}
We can now claim that the relation (\ref{eq: controlbornee}) given in proposition \ref{estimclass} holds true, using that
$$ e^{\Gamma C_{t}} |U_{t}^{n}|^{2} \le \mb{E}\left(\disp{\sup_{ t\le u\le T}e^{\Gamma C_{u}} |U_{u}^{n}|^{2}|\mc{F}_{t} }\right).$$
To prove the boundedness in $S^{\infty}$, we rely on the two following properties: on the one hand, $|g^{n}(u, 0, 0)| \leq \bar{\alpha}_{u}$, with the process $\bar{\alpha}$ satisfying:  $\disp{\int_{0}^{T}\bar{\alpha}_{s}dC_{s} \le a < \infty}$, $\mb{P}$-a.s. and, on the other hand and for all $n$, the random variable $U_{T}^{n} :=e^{\beta B}$ is bounded.\\
\begin{flushright}
$\square$
\end{flushright}


\begin{thebibliography}{}
\bibitem{Stricker2}
Ansel, J.-P. and Stricker, C.,
 \newblock {\em Lois de martingale, densit\'es et d\'ecomposition de
              {F}\"ollmer-{S}chweizer},
\newblock {\em Ann. Inst. H. Poincar\'e Probab. Statist.}, \textbf{28}(3) :375--392, 1992.

\bibitem{Becherer}
Becherer, D.,
\newblock{\em Bounded solutions to Backward SDE's with jumps for utility optimization and indifference hedging},
\newblock{ \em Ann. Appl. Probab.}, \textbf{16}(4) : 2027--2054, 2006.

\bibitem{Fritelli}
Biagini, S. and Frittelli, M.,
  \newblock {\em Utility maximization in incomplete markets for unbounded
              processes},
  \newblock {\em  Finance Stoch.}, \textbf{9}(4) : 493--517, 2005.


\bibitem{BriandetCoquet}
Briand, P., Coquet, F., Hu, Y., M{\'e}min, J. and Peng, S.,
 \newblock {\em A converse comparison theorem for {BSDE}s and related
              properties of {$g$}-expectation},
\newblock {\em    Electron. Comm. Probab.}, \textbf{5} : 101--117, 2000.
 
   
\bibitem{BriandetHu}
 Briand, P. and Hu, Y.,
 \newblock {\em BSDE with quadratic growth and unbounded terminal value},
 \newblock {\em Probab. Theory Related Fields}, \textbf{136}(4) : 604--618, 2006.

\bibitem{DelbaenetSch}
Delbaen, F. and Schachermayer, W.,
 \newblock {\em The existence of absolutely continuous local martingale
              measures},
 \newblock {\em  Ann. Appl. Probab.}, \textbf{5}(4) : 926--945, 1995.

 \bibitem{ElkaretMazliak}
El Karoui, N. and Huang, S.-J.,
\newblock {\em A general result of existence and uniqueness of backward
              stochastic differential equations},
\newblock{\em Backward stochastic differential equations, Pitman Res. Notes Math. Ser.}, \textbf{364} : 27--36, Longman, Harlow, 1997.


\bibitem{ElKaretPeng}
El Karoui, N., Peng S. and Quenez M.C.,
\newblock {\em Backward stochastic differential equations in finance},
\newblock{\em Math. Finance}, \textbf{7}(1) : 1--71, 1997.

\bibitem{ElKaretRouge}
El Karoui, N. and Rouge, R.,
\newblock {\em Pricing via utility maximization and entropy},
\newblock{\em Math. Finance}, \textbf{10}(2) : 259--276, 2000. 
 
 \bibitem{HPSCHW01}
Heath, D. and Platen, E. and Schweizer, M.,
 \newblock {\em A comparison of two quadratic approaches to hedging in
              incomplete markets},
  \newblock {\em Math. Finance}, \textbf{11}(4): 385--413, 2001.


\bibitem{ImkelleretHu}
Hu, Y., Imkeller, P. and M\"uller, M.,
\newblock{\em Utility maximization in incomplete markets},
\newblock{\em Ann. Appl. Probab.}, \textbf{15}(3) : 1691--1712, 2005.

\bibitem{Kazamaki}
Kazamaki, N.,
 \newblock {\em Continuous Exponential Martingales and {BMO}},
\newblock{\em   Lecture Notes in Math.}, \textbf{1579}, Springer, Berlin, 1994. 
 
    
\bibitem{mkobylanski}
Kobylanski, M.,
\newblock {\em Backward stochastic differential equations and partial
              differential equations with quadratic growth},
\newblock {\em Ann. Probab.}, \textbf{28}(2) : 558--602, 2000.
  
\bibitem{LepeletSanm}
Lepeltier, J. P. and San Martin, J.,
\newblock {\em Existence for {BSDE} with superlinear-quadratic coefficient},
\newblock {\em Stochastics Stochastics Rep.}, \textbf{63}(3-4) : 227--240, 1998.

\bibitem{ManiaetSchw}
Mania, M. and Schweizer, M.,
  \newblock {\em    Dynamic exponential utility indifference valuation},
 \newblock {\em  Ann. Appl. Probab.}, \textbf{15}(3) : 2113--2143, 2005.

\bibitem{Owen}
Owen, M. P.,
 \newblock{ \em Utility based optimal hedging in incomplete markets},
  \newblock {\em Ann. Appl. Probab.},
  \textbf{12}(2) : 691--709 , 2002.

\bibitem{Pardoux2}  
Pardoux, {\'E}.,
 \newblock {\em  B{SDE}s, weak convergence and homogenization of semilinear
              {PDE}s},
 \newblock {\em Nonlinear analysis, differential equations and control}, \textbf{528} : 503--549, 1999.


\bibitem{ParetPeng}
 Pardoux, {\'E}. and Peng, S.,
 \newblock {\em Adapted solution of a backward stochastic differential
              equation},
 \newblock {\em Systems Control Lett.} \textbf{14}(1) : 55--61, 1990.
\bibitem{Protter}
Protter, P.,
 \newblock {\em Stochastic Integration and Differential Equations}, Springer, Berlin, 2004.
 
 \bibitem{RevuzYor}
Revuz, D. and Yor, M.,
\newblock {\em Continuous Martingales and {B}rownian Motion}, Springer, Berlin, 1999.


\bibitem{Schach1}
Schachermayer, W.,
\newblock {\em Utility maximization in incomplete markets},
\newblock {\em Lecture Notes in Math.}, \textbf{1856} : 255--293, Springer, Berlin, 2004.
   
\bibitem{Zhariphopoulou}
Zhariphopoulou, T.,
\newblock{ \em A solution approach to valuation with unhedgeable risk}
\newblock{\em Finance and Stochastics}, \textbf{5} : 61--82, Springer, 2001.
 \end{thebibliography}
\end{document}